\documentstyle{article}

\def\aA{\mbox{$\mathcal A$}}
\def\D{\mbox{$\mathcal D$}}
\def\cirk{\,{\raisebox{.3ex}{\tiny $\circ$}}\,}
\def\kon{\wedge}
\def\str{\rightarrow}
\def\rts{\leftarrow}
\def\mj{\mbox{\bf 1}}
\def\koc{{\raisebox{-.2ex}{$\Box$}}}
\def\pl{\!+\!}
\def\mn{\!-\!}
\def\prop#1#2{\vspace{2ex} \noindent{\sc #1.} {\it #2} \par \vspace{2ex}}
\def\dkz{\noindent{\sc Proof. }}
\def\qed{\hfill $\dashv$\vspace{2ex}}
\def\Cm{\mbox{$\mathbf{C^m}$}}

\begin{document}

\title{Medial Commutativity}
\author{\small {\sc Kosta Do\v sen} and {\sc Zoran Petri\' c}
\\[1ex]
{\small Mathematical Institute, SANU}\\[-.5ex]
{\small Knez Mihailova 35, p.f. 367, 11001 Belgrade,
Serbia}\\[-.5ex]
{\small email: \{kosta, zpetric\}@mi.sanu.ac.yu}}
\date{}
\maketitle

\vspace{-6ex}

\begin{abstract}
\noindent It is shown that all the assumptions for symmetric
monoidal categories flow out of a unifying principle involving
natural isomorphisms of the type ${(A\kon B)\kon(C\kon
D)\str(A\kon C)\kon(B\kon D)}$, called medial commutativity.
Medial commutativity in the presence of the unit object enables us
to define associativity and commutativity natural isomorphisms. In
particular, Mac Lane's pentagonal and hexagonal coherence
conditions for associativity and commutativity are derived from
the preservation up to a natural isomorphism of medial
commutativity by the biendofunctor $\kon$. This preservation boils
down to an isomorphic representation of the Yang-Baxter equation
of symmetric and braid groups. The assumptions of monoidal
categories, and in particular Mac Lane's pentagonal coherence
condition, are explained in the absence of commutativity, and also
of the unit object, by a similar preservation of associativity by
the biendofunctor $\kon$. In the final section one finds coherence
conditions for medial commutativity in the absence of the unit
object. These conditions are obtained by taking the direct product
of the symmetric groups $S_{n \choose i}$ for $0\leq i\leq n$.

\end{abstract}

\vspace{.3cm}

\noindent {\small {\it Mathematics Subject Classification} ({\it
2000}): 18D10, 19D23, 20B30}

\vspace{.5ex}

\noindent {\small {\it Keywords$\,$}: associativity,
commutativity, monoidal categories, symmetric monoidal categories,
coherence, Mac Lane's pentagon, Mac Lane's hexagon, Yang-Baxter
equation, symmetric groups, binomial coefficients}

\vspace{.5ex}

\noindent {\small {\it Acknowledgement$\,$}. Work on this paper
was supported by the Ministry of Science of Serbia (Grant
144013).}

\vspace{1cm}

\vspace{-7ex}

\section{\large\bf Introduction}
The purpose of the main section of this paper (Section~2) is to
show that the notion of symmetric monoidal category may be so
presented that all the assumptions for it flow out of a unifying
principle. Among these assumptions we find in particular that
$\kon$ (often written $\otimes$) is a biendofunctor of the
category \aA, i.e.\ a functor from $\aA^2$ (which is the product
category ${\aA\times\aA}$) to \aA, that preserves up to a natural
isomorphism the structure of $\aA^2$ together with the induced
biendofunctor $\kon^2$ of $\aA^2$ defined coordinatewise in terms
of $\kon$. This is related to the notion of monoidal functor (see
\cite{EK66}, Sections II.1 and III.1, \cite{JS93}, \cite{ML71},
second edition, Section XI.2, and \cite{DP04}, Section 2.8). The
natural isomorphism $c^m$ required for that has as components the
arrows of \aA
\[
c^m_{A,B,C,D}:(A\kon B)\kon(C\kon D)\str(A\kon C)\kon(B\kon D).
\]

We call \emph{medial commutativity} the principle on which $c^m$
is based. In universal algebra, besides being called
\emph{medial}, it was called \emph{abelian}, \emph{alternation},
\emph{bicommutative}, \emph{bisymmetric}, \emph{entropic},
\emph{surcommutative} and \emph{transposition} (see \cite{JK83};
early references for medial commutativity in category theory are
\cite{EH62}, Section~3, and \cite{EK66}, Section III.3, where it
was called \emph{middle-four interchange}---see also \cite{DP96},
Section~2, and \cite{DP04}, Section 9.1; medial commutativity is a
kind of riffle shuffle). In the presence of the unit object $\top$
(often written $\mathrm I$), one defines in terms of $c^m$ the
natural isomorphisms whose components are of the types
${A\kon(C\kon D)\str(A\kon C)\kon D}$ and ${B\kon C\str C\kon B}$,
i.e.\ associativity and commutativity.

Next we introduce a general notion of what it means for a functor
to preserve a natural transformation, and we require that $\kon$
preserve $c^m$, the natural isomorphism it has previously induced.
Together with similar assumptions involving $\top$, this delivers
exactly the notion of symmetric monoidal category. In particular,
Mac Lane's pentagonal and hexagonal coherence conditions for
associativity and commutativity are derived from the preservation
of $c^m$ by $\kon$ (in the presence of $\top$).

The preservation of $c^m$ by $\kon$ boils down to an isomorphic
representation of the Yang-Baxter equation of symmetric and braid
groups. This may be taken as an explanation of the Yang-Baxter
equation through our notion of preservation. We deal with this
matter briefly in Section~3. It would lead us to far to go into
the details of this representation, which stems from \cite{B37},
and we give only summary indications, pointing to references that
cover the subject more exhaustively. We also mention at the end of
that section a parallel that exists between matters treated here
and a result of \cite{DP06}.

In Section~4 of the paper we explain in the same spirit the
assumptions for monoidal categories, with and without $\top$. We
show in particular how in the absence of symmetry, i.e.\
commutativity, and $\top$, Mac Lane's pentagonal coherence
condition for associativity may be derived from a kind of
preservation of associativity by $\kon$.

In the final section of the paper (Section~5) we show what
assumptions for $c^m$ in the absence of $\top$ are necessary and
sufficient for coherence. The equations required are obtained by
taking the direct product of the symmetric groups $S_{n \choose
i}$ for $0\leq i\leq n$.

We presuppose for this paper an acquaintance with the notions of
monoidal and symmetric monoidal category and with Mac Lane's
coherence results for these categories (see \cite{ML63},
\cite{ML71}, or \cite{DP04}). We also rely on some standard
notions of category theory, whose definitions may be found in
\cite{ML71}. A modest acquaintance with symmetric groups, and in
particular with their standard presentation by transpositions of
immediate neighbours, is required for Section~5. For motivation we
also rely, but not excessively, on some notions of \cite{DP04}
(Section 2.8).

We understand coherence in the sense of Mac Lane's result for
symmetric monoidal categories. His theorem may be understood as
asserting the existence of a faithful functor from the symmetric
monoidal category freely generated by a set of objects to the
category whose arrows are permutations of finite ordinals (see
\cite{DP04}, Sections 1.1 and 2.9, for a general approach to
coherence in this spirit).

To fix terminology, we call \emph{categorial} equations the
following usual equations assumed for categories:
\begin{tabbing}
\mbox{\hspace{2.5em}}\=$({\mbox{{\it
cat}~1}})$\mbox{\hspace{2em}}\=$f\cirk \mj_A=\mj_B\cirk f=f,
\quad{\mbox{\rm for}}\; f\!:A\str B,$
\\*[1ex]
\>$({\mbox{{\it cat}~2}})$\>$h\cirk (g\cirk f)=(h\cirk g)\cirk f.$
\end{tabbing}
We call \emph{bifunctorial} equations for $\kon$ the equations
\begin{tabbing}
\mbox{\hspace{2.5em}}\=$({\mbox{{\it
cat}~1}})$\mbox{\hspace{2em}}\=$f\cirk \mj_A=\mj_B\cirk f=f,
\;{\mbox{\rm for}}\; f\!:A\str B,$\kill

\>$({\mbox{{\it bif}~1}})$\>$\mj_A\kon\mj_B=\mj_{A\kon B},$
\\*[1ex]
\>$({\mbox{{\it bif}~2}})$\>$(g_1\cirk f_1)\kon (g_2\cirk
f_2)=(g_1\kon g_2)\cirk (f_1\kon f_2).$
\end{tabbing}
The \emph{naturality} equation for $c^m$ is
\begin{tabbing}
\mbox{\hspace{2.5em}}\=$({\mbox{{\it
cat}~1}})$\mbox{\hspace{2em}}\=$f\cirk \mj_A=\mj_B\cirk f=f,
\;{\mbox{\rm for}}\; f\!:A\str B,$\kill

\>$(c^m\;{\mbox\it nat})$\>$((f\kon h)\kon(g\kon j))\cirk
c^m_{A,B,C,D}=c^m_{A',B',C',D'}\cirk((f\kon g)\kon(h\kon j)).$
\end{tabbing}
The \emph{isomorphism} equation for $c^m$ is the equation
$(c^mc^m)$ given in the next section. We have analogous naturality
and isomorphism equations for other natural isomorphisms to be
encountered in the text.

\section{\large\bf Symmetric monoidal categories}
For $m\geq 1$, an $m$-endofunctor of a category \aA\ is a functor
from the product category $\aA^m$ to \aA. As in the introduction
above, we call 2-endofunctors \emph{biendofunctors}, and
1-endofunctors just \emph{endofunctors}. We can take a special
object of \aA\ to be a 0-endofunctor of \aA, taking that $\aA^0$
is a trivial one-object category.

We will make in this section a number of suppositions that will
make out of \aA\ at the end a symmetric monoidal category. The
first supposition is the following:
\begin{quote}
(I)\quad $\kon$ is a biendofunctor of the category \aA.
\end{quote}
We write as usual $a\kon b$, rather than $\kon(a,b)$, for
$a,b,\ldots$ being either objects or arrows. We reserve
$A,B,\ldots$ for objects, and $f,g,\ldots$ for arrows.

We define in terms of $\kon$ a biendofunctor $\kon^2$ of $\aA^2$
by stipulating, in a coordinatewise manner, that
\[
(a_1,a_2)\kon^2(b_1,b_2)=_{df}(a_1\kon b_1,a_2\kon b_2).
\]
Note that the straightforward definition
\[
(a_1,a_2){\kon^2}'(b_1,b_2)=_{df}(a_1\kon a_2,b_1\kon b_2)
\]
would also define a biendofunctor of $\aA^2$, but $\kon^2$ has
properties that the straightforwardly defined ${\kon^2}'$ does not
have. One of these properties is the following.

Take the \emph{diagonal} subcategory \D\ of $\aA^2$, which is the
image in $\aA^2$ of the diagonal functor $D$ from \aA\ to $\aA^2$,
defined by ${Da=_{df}(a,a)}$. Then the coordinatewise definition
of $\kon^2$ gives a biendofunctor of \D, while the straightforward
definition of ${\kon^2}'$ does not deliver that. The category \D\
is a category isomorphic to \aA\ such that ${Da\kon^2 Db= D(a\kon
b)}$. (Another property that distinguishes $\kon^2$ from
${\kon^2}'$ will appear below with the functor $(M,M)$ and the
natural transformation $(\alpha,\alpha)$.)

Suppose, for $i\in\{1,2\}$, that $M_i$ is an $m$-endofunctor of
the category $\aA_i$. We say that a functor $F$ from $\aA_1$ to
$\aA_2$ is an \emph{upward functor} from $(\aA_1,M_1)$ to
$(\aA_2,M_2)$ when there is a natural transformation $\psi$ whose
components are the arrows of $\aA_2$
\[
\psi_{A_1,\ldots,A_m}\!:M_2(FA_1,\ldots,FA_m)\str
FM_1(A_1,\ldots,A_m).
\]

We assume, of course, that $\psi$ is natural in the indices
$A_1,\ldots,A_m$. Let the functor ${F^m}'$ from $\aA_1^m$ to
$\aA_2^m$ be defined out of $F$ in the straightforward manner;
namely,
\[
{{F^m}'(a_1,\ldots,a_m)=_{df}(Fa_1,\ldots,Fa_m)}
\]
(cf.\ above the definition of ${\kon^2}'$). Then $\psi$ is a
natural transformation from the composite functor $M_2{F^m}'$ to
the composite functor $FM_1$, which go from $\aA_1^m$ to $\aA_2$.
It should be clear from the context from what functor to what
functor a natural transformations such as $\psi$ goes, and we will
not always mention such things explicitly.

We say for an upward functor that it is \emph{loyal} when $\psi$
is a natural isomorphism. When this natural isomorphism is
identity, then $F$ preserves $M_1$ exactly, as one would expect in
adapting the notion homomorphism to categories..

The second supposition we make in order to make out of \aA\ a
symmetric monoidal category is the following:
\begin{quote}
(II)\quad $\kon$ is an upward functor from $(\aA^2,\kon^2)$ to
$(\aA,\kon)$.
\end{quote}
The natural transformation $\psi$ of this upward functor is $c^m$,
whose components are the arrows of \aA
\[
c^m_{A_1,A_2,A_3,A_4}\!:(A_1\kon A_2)\kon(A_3\kon A_4)\str(A_1\kon
A_3)\kon(A_2\kon A_4).
\]

The source of $c^m_{A_1,A_2,A_3,A_4}$ is
${\kon(\kon(A_1,A_2),\kon(A_3,A_4))}$, and the target is
${\kon\kon^2((A_1,A_2),(A_3,A_4))}$. The natural transformation
$c^m$ goes from the composite functor $\kon{\kon^2}'$ to the
composite functor $\kon\kon^2$; here $M_2$ and $F$ are both $\kon$
and $M_1$ is $\kon^2$.

Our third supposition is the following:
\begin{quote}
(III)\quad \aA\ has a special object $\top$.
\end{quote}
As we said at the beginning of the section, we may take this
special object as a 0-endofunctor of \aA.

In terms of the biendofunctor $\kon$ and of $\top$ we may define
the endofunctor $\kon\top$ of \aA\ by $\kon\top(a)=_{df}
a\kon\top$, where if $a$ is an arrow, then $\top$ stands for
$\mj_\top$.

For the fourth supposition we take $I$ to be the identity functor
of \aA:
\begin{quote}
(IV$\delta$)\quad $I$ is a loyal functor from $(\aA,\kon\top)$ to
$(\aA,I)$.
\end{quote}
The natural isomorphism $\psi$ of this loyal functor is
$\delta^\str$, whose components are the arrows of \aA
\[
\delta^\str_A\!:A\kon\top\str A, \quad {\mbox{\rm with
inverses}}\quad \delta^\rts_A\!:A\str A\kon\top.
\]

In terms of the biendofunctor $\kon$ and of $\top$ we may also
define the endofunctor $\top\kon$ of \aA\ by $\top\kon(a)=_{df}
\top\kon a$, and the second part of the fourth supposition is:
\begin{quote}
(IV$\sigma$)\quad $I$ is a loyal functor from $(\aA,\top\kon)$ to
$(\aA,I)$.
\end{quote}
The natural isomorphism $\psi$ of this loyal functor is
$\sigma^\str$, whose components are the arrows of \aA
\[
\sigma^\str_A\!:\top\kon A\str A, \quad {\mbox{\rm with
inverses}}\quad \sigma^\rts_A\!:A\str \top\kon A.
\]

Before we can formulate the next supposition, we must introduce a
number of notions. We define inductively the notion of
\emph{shape} and of its arity in the following manner:
\begin{itemize}
\item[] $\top$ is a shape of arity 0, and \koc\ is a shape of
arity 1; \item[] if $M$ and $N$ are shapes of arities $m$ and $n$
respectively, then $(M\kon N)$ is a shape of arity $m\pl n$.
\end{itemize}
We take the outermost parentheses of shapes for granted, and omit
them.

A shape $M$ of arity $m$ defines an $m$-endofunctor of \aA, such
that $M(a_1,\ldots,a_m)$ is obtained by putting $a_i$, where
$1\leq i\leq m$, for the $i$-th $\koc$, counting from the left, in
the shape $M$. (For arrows, we read $\top$ as $\mj_\top$, as we
said above.) An $m$-endofunctor defined by a shape $M$ and a
permutation $\pi$ of $\{1,\ldots,m\}$ define an $m$-endofunctor
$M^\pi$ such that
\[
M^\pi(a_1,\ldots,a_m)=_{df}M(a_{\pi(1)},\ldots,a_{\pi(m)}).
\]

We shall use the following abbreviations:
\begin{tabbing}
\mbox{\hspace{8em}}\=$\vec{a}_m$ \mbox{\hspace{3em}}\= for \quad
${a_1,\ldots,a_m}$,\\[1ex]
\>$\vec{a}_{\pi(m)}$ \> for \quad
${a_{\pi(1)},\ldots,a_{\pi(m)}}$,\\[1ex]
\>$\overrightarrow{a_m, a_m'\,}$ \> for \quad ${a_1,
a_1',\ldots,a_m, a_m'}$,\\[1ex]
\>$\overrightarrow{(a_m, a_m')}$ \> for \quad ${(a_1,
a_1'),\ldots,(a_m, a_m')}$,\\[1ex]
\>$\overrightarrow{a_m\kon a_m'\,}$ \> for \quad ${a_1\kon
a_1',\ldots,a_m\kon a_m'}$,
\end{tabbing}
and other analogous abbreviations, made on the same pattern.

Next we give by induction on the complexity of the shape $M$ of
arity $m$ the definition of the natural transformation $\psi^M$,
whose components are the arrows of \aA
\[
\psi^M_{\scriptsize\overrightarrow{A_m,A_m'}}\!\!:M(\overrightarrow{A_m\kon
A_m'\,})\str M(\vec{A}_m)\kon M(\vec{A}_m').
\]
This is a natural transformation from the $2m$-endofunctor of \aA\
defined by the shape $M(\koc\kon\koc,\ldots,\koc\kon\koc)$ to the
$2m$-endofunctor of \aA\ defined by the shape
$M(\koc,\ldots,\koc)\kon M(\koc,\ldots,\koc)$ and a permutation of
$\{1,\ldots,2m\}$. Here is the inductive definition:
\begin{tabbing}
\mbox{\hspace{2em}}\=$\psi^\top=\delta^\rts_\top:\top\str\top\kon\top$,\\[1.5ex]
\>$\psi^\Box_{A,A'}=\mj_{A\kon A'}$,\\[2ex]
\>$\psi^{M\kon
N}_{\scriptsize\overrightarrow{A_m,A_m'},\overrightarrow{B_n,B_n'\,}}=
c^m_{M(\vec{A}_m), M(\vec{A}_m'),N(\vec{B}_n),N(\vec{B}_n')}\cirk
(\psi^M_{\scriptsize\overrightarrow{A_m,A_m'}}\kon
\psi^N_{\scriptsize\overrightarrow{B_n,B_n'\,}})$.
\end{tabbing}
In the last clause, $M$ is a shape of arity $m$ and $N$ a shape of
arity $n$.

If in the shape $M$ there are no occurrences of $\top$, then
$\delta^\rts$ is absent from $\psi^M$. Note that
\[
\psi^{\Box\kon\Box}_{A,A',B,B'}=c^m_{A,A',B,B'}.
\]

Let $\kon^m$, for $m\geq 1$, be the biendofunctor of $\aA^m$
defined coordinatewise by
\[
\vec{a}_m\kon^m\vec{a'}_m=_{df}(\overrightarrow{a_m\kon a_m'}).
\]
This generalizes the definition of $\kon^2$ given previously, and
$\kon$ is $\kon^1$. Let $(M,M)$ be the functor from $(\aA^m)^2$ to
$\aA^2$ defined straightforwardly by
\[
(M,M)((\vec{a}_m),(\vec{a'}_m))=_{df}(M(\vec{a}_m),M(\vec{a'}_m)).
\]
Then $\psi^M$ may be conceived as a natural transformation between
two functors from $(\aA^m)^2$ to \aA. The source functor here is
the composite functor $M\kon^m$, and the target functor is the
composite functor $\kon(M,M)$. In that case the components of
$\psi^M$ should be indexed as follows:
\[
\psi^M_{(\vec{A}_m),(\vec{A}_m')}.
\]

Another possibility is to take the functor ${\kon^m}'$ from
$(\aA^2)^m$ to $\aA^m$, whose straightforward definition
\[
{\kon^m}'(\overrightarrow{(a_m,a_m')})=_{df}
(\overrightarrow{a_m\kon a_m'})
\]
has the same right-hand side as the definition of $\kon^m$ above.
Let $M^2$ be obtained from the shape $M$ by substituting $\kon^2$
for $\kon$ and $(\top,\top)$ for $\top$. Note that for the
$m$-endofunctor of $\aA^2$ defined by $M^2$, according to the
coordinatewise definition of $\kon^2$, we have
\[
M^2(\overrightarrow{(a_m,a_m')})= (M(\vec{a}_m),M(\vec{a'}_m)),
\]
which has the same right-hand side as the definition of $(M,M)$
above. Then $\psi^M$ may be conceived as a natural transformation
between two functors from $(\aA^2)^m$ to \aA. The source functor
here is the composite functor $M{\kon^m}'$, and the target functor
is the composite functor $\kon M^2$. In that case the components
of $\psi^M$ should be indexed as follows:
\[
\psi^M_{\scriptsize\overrightarrow{(A_m,A_m')}},
\]
and our official indexation is obtained from this one just by
deleting the parentheses. This indexation is in tune with the
source of the components of $\psi^M$.

Let ${M^2}'$ be obtained from $M^2$ by substituting ${\kon^2}'$
for $\kon^2$ (see the beginning of the section for the definition
of ${\kon^2}'$). Then ${M^2}'$ defines an $m$-endofunctor of
$\aA^2$ that does not amount to $(M,M)$, as $M^2$ does. For
example,
\[
((\koc{\kon^2}'\koc){\kon^2}'\koc)((a_1,a_1'),(a_2,a_2'),(a_3,a_3'))=
((a_1\kon a_1')\kon(a_2\kon a_2'), a_3\kon a_3').
\]

Let $M_i$, for $i\in\{1,2\}$, be shapes of arity $m$; so they
define $m$-endofunctors of \aA. Let $\alpha$ be a natural
transformation from $M_1$ to $M_2$. Then there exists a natural
transformation $(\alpha,\alpha)$ from the functor $(M_1,M_1)$ to
the functor $(M_2,M_2)$. The components of $(\alpha,\alpha)$ are
pairs of components of $\alpha$. As we saw above, the natural
transformation $(\alpha,\alpha)$ may be conceived as going from
$M_1^2$ to $M_2^2$. This does not work if we take ${M_i^2}'$
instead of $M_i^2$.

Let now $\alpha$ be a natural transformation from the
$m$-endofunctor of \aA\ defined by the shape $M_1$ to the
$m$-endofunctor of \aA\ defined by the shape $M_2$ and a
permutation $\pi$ of $\{1,\ldots,m\}$. We say that $\alpha$ is
\emph{upward preserved} by $\kon$ when diagrams of the following
form commute in \aA:
\begin{center}
\begin{picture}(300,105)
\put(3,10){$M_1(\vec{A}_m)\kon M_1(\vec{A}_m')$}

\put(8,80){$M_1(\overrightarrow{A_m\kon A_m'})$}

\put(184,80){$M_2(\overrightarrow{A_{\pi(m)}\kon A_{\pi(m)}'})$}

\put(180,10){$M_2(\vec{A}_{\pi(m)})\kon M_2(\vec{A}_{\pi(m)}')$}

\put(7,48){$\psi^{M_1}_{\scriptsize\overrightarrow{A_m,A_m'}}$}

\put(237,48){$\psi^{M_2}_{\scriptsize\overrightarrow{A_{\pi(m)},A_{\pi(m)}'}}$}

\put(108,6){$\alpha_{\vec{A}_m}\kon\alpha_{\vec{A}_m'}$}

\put(92,15){\vector(1,0){80}}

\put(110,93){$\alpha_{\scriptsize\overrightarrow{A_m\kon A_m'}}$}

\put(92,85){\vector(1,0){80}}

\put(46,70){\vector(0,-1){40}}

\put(233,70){\vector(0,-1){40}}

\end{picture}
\end{center}
i.e., we have in \aA\ the equation
\begin{tabbing}
\mbox{\hspace{2.5em}}$(\psi\alpha)\quad\quad
\psi^{M_2}_{\scriptsize\overrightarrow{A_{\pi(m)},A_{\pi(m)}'}}\cirk
\alpha_{\scriptsize\overrightarrow{A_m\kon
A_m'}}=(\alpha_{\vec{A}_m}\kon\alpha_{\vec{A}_m'})\cirk
\psi^{M_1}_{\scriptsize\overrightarrow{A_m,A_m'}}$.
\end{tabbing}
(This equation is an instance of the equation $(\psi\alpha)$ of
\cite{DP04}, Section 2.8; something analogous in a concrete
context, different from ours, occurs in \cite{BFSV03}, Section~1.)

We can now formulate our next supposition in three parts. We have
first:
\begin{quote}
(V$c^m$)\quad $c^m$ is upward preserved by $\kon$.
\end{quote}
This means that we have in \aA\ the following instance of
$(\psi\alpha)$:
\begin{tabbing}
\mbox{\hspace{2.5em}}$(\psi c^m)$\quad\= $c^m_{A_1\kon
A_3,A_1'\kon A_3',A_2\kon A_4,A_2'\kon
A_4'}\cirk(c^m_{A_1,A_1',A_3,A_3'}\kon
c^m_{A_2,A_2',A_4,A_4'})\cirk$\\*[.5ex] \`$\cirk c^m_{A_1\kon
A_1',A_2\kon A_2',A_3\kon A_3',A_4\kon A_4'}=$\\*[1ex]
\>$(c^m_{A_1,A_2,A_3,A_4}\kon c^m_{A_1',A_2',A_3',A_4'})\cirk
c^m_{A_1\kon A_2,A_1'\kon A_2',A_3\kon A_4,A_3'\kon
A_4'}\cirk$\\*[.5ex] \`$\cirk(c^m_{A_1,A_1',A_2,A_2'}\kon
c^m_{A_3,A_3',A_4,A_4'})$
\end{tabbing}
(see the next section for a graphical presentation of this
equation, which is obtained from the hexagonal interchange
equation of \cite{BFSV03}, end of Definition 1.7, by deleting the
indices $i$, $j$ and $k$).

The second part of the fifth supposition is:
\begin{quote}
(V$\delta$)\quad $\delta^\str$ is upward preserved by $\kon$.
\end{quote}
This means that we have in \aA\ the following equation:
\begin{tabbing}
\mbox{\hspace{2.5em}}$(\psi\delta)\quad$\= $\delta^\str_{A\kon
A'}=(\delta^\str_A\kon\delta^\str_{A'})\cirk
c^m_{A,A',\top,\top}\cirk(\mj_{A\kon A'}\kon\delta^\rts_\top)$.
\end{tabbing}

The final part of the fifth supposition, analogous to (V$\delta$),
is:
\begin{quote}
(V$\sigma$)\quad $\sigma^\str$ is upward preserved by $\kon$.
\end{quote}
This means that we have in \aA\ the following equation:
\begin{tabbing}
\mbox{\hspace{2.5em}}$(\psi\sigma)\quad$\= $\sigma^\str_{A\kon
A'}=(\sigma^\str_A\kon\sigma^\str_{A'})\cirk
c^m_{\top,\top,A,A'}\cirk(\delta^\rts_\top\kon\mj_{A\kon A'})$.
\end{tabbing}

As a consequence of (V$\delta$) and (V$\sigma$) we have that
$\delta^\rts$ and $\sigma^\rts$ are upward preserved by $\kon$. It
is clear that \mj\ understood as a natural transformation from the
identity functor to the identity functor is upward preserved by
$\kon$. Then it can be shown by induction on the complexity of $M$
that for every $M$ the natural transformation $\psi^M$ is upward
preserved by $\kon$. As a matter of fact, something even more
general holds. Every natural transformation defined in terms of
\mj, $c^m$, $\delta^\str$, $\delta^\rts$, $\sigma^\str$,
$\sigma^\rts$, $\kon$ and $\cirk$ is upward preserved by $\kon$.

As a consequence of $(\psi\delta)$ we have
\begin{tabbing}
\mbox{\hspace{2.5em}}$(c^m34)$\mbox{\hspace{2em}}\=
$c^m_{A_1,A_2,\top,\top}=(\delta^\rts_{A_1}\kon\delta^\rts_{A_2})\cirk
\delta^\str_{A_1\kon A_2}\cirk(\mj_{A_1\kon
A_2}\kon\delta^\str_\top)$
\end{tabbing}
---an equation we will need later for proving that \aA\ is a
symmetric monoidal category. We also have
\begin{tabbing}
\mbox{\hspace{2.5em}}\=$(c^m34)$\mbox{\hspace{2em}}\=
$c^m_{A_1,A_2,\top,\top}$\kill \>$(c^m234)$\>
$c^m_{A,\top,\top,\top}=\mj_{(A\kon\top)\kon(\top\kon\top)}$,
\end{tabbing}
for which we use $(c^m34)$ and the following instance of the
naturality equation for $\delta^\str$:
\[
\delta^\str_A\cirk\delta^\str_{A\kon\top}=\delta^\str_A\cirk(\delta^\str_A
\kon\mj_\top),
\]
which delivers $\delta^\str_{A\kon\top}=\delta^\str_A
\kon\mj_\top$.

Our penultimate supposition is:
\begin{quote}
(VI)\quad $c^m$ is an isomorphism inverse to itself.
\end{quote}
This means that we have in \aA\ the following equation:
\begin{tabbing}
\mbox{\hspace{2.5em}}\=$({\mbox{{\it
cat}~1}})$\mbox{\hspace{2em}}\=$f\cirk \mj_A=\mj_B\cirk f=f,
\;{\mbox{\rm for}}\; f\!:A\str B,$\kill

\>$(c^m c^m)$\>$c^m_{A_1,A_3,A_2,A_4}\cirk
c^m_{A_1,A_2,A_3,A_4}=\mj_{(A_1\kon A_2)\kon(A_3\kon A_4)}.$
\end{tabbing}
Together with supposition (II) this implies that $\kon$ is a loyal
functor from $(\aA^2,\kon^2)$ to $(\aA,\kon)$, but it amounts to
more than just supposing loyalty. We do not suppose only that
$c^m$ has an inverse, but also that this inverse is $c^m$ itself.
Supposing that $c^m$ is an isomorphism without supposing that it
is inverse to itself points towards braided monoidal categories
(see \cite{JS93} and \cite{ML71}, second edition). For the time
being however we leave aside the investigation of braided monoidal
categories in the spirit of this paper, and deal in the paper only
with the simpler symmetric monoidal categories.

From $(c^m34)$ and $(c^m c^m)$ we derive
\begin{tabbing}
\mbox{\hspace{2.5em}}$(c^m24)$\mbox{\hspace{2em}}\=
$c^m_{A_1,\top,A_3,\top}=(\mj_{A_1\kon
A_3}\kon\delta^\rts_\top)\cirk \delta^\rts_{A_1\kon
A_3}\cirk(\delta^\str_{A_1}\kon\delta^\str_{A_3})$.
\end{tabbing}

The last supposition we make in order to make out of \aA\ a
symmetric monoidal category is the following:
\begin{quote}
(VII)\mbox{\hspace{2em}} $\delta^\str_\top=\sigma^\str_\top$.
\end{quote}
This supposition, which is equivalent to
${\delta^\rts_\top=\sigma^\rts_\top}$, may be justified as
follows. In defining inductively $\psi^M$ above we have defined
$\psi^\top$ as ${\delta^\rts_\top\!:\top\str\top\kon\top}$, but we
could equally well have chosen to define $\psi^\top$ in a
different manner, as ${\sigma^\rts_\top\!:\top\str\top\kon\top}$.
With (VII) we ensure that there is no difference between these
definitions---we abolish the arbitrary favouring of
$\delta^\rts_\top$ over $\sigma^\rts_\top$.

An alternative, longer, way to justify (VII) is to require that
$\top$, conceived as nullary $\kon$, is an upward functor from
$(\aA^0,\kon)$ to $(\aA,\kon)$, where for $\ast$ being the unique
object of $\aA^0$ we have $\ast\kon\ast=\ast$. This requires that
there be in \aA\ an arrow
$\psi_{\ast,\ast}\!:\top\kon\top\str\top$. For every shape $M$ we
define by induction an arrow
$\psi^M\!:M(\top,\ldots,\top)\str\top$:
\begin{tabbing}
\mbox{\hspace{2.5em}}$(c^m24)$\mbox{\hspace{2em}}\=\kill
\>$\psi^\top=\psi^\Box=\mj_\top$,\\[1ex]
\>$\psi^{M\kon N}=\psi_{\ast,\ast}\cirk(\psi^M\kon\psi^N)$.
\end{tabbing}
For $\alpha$ being as for $(\psi\alpha)$, let us say that $\alpha$
is \emph{upward preserved} by $\top$ when in \aA\ we have the
equation
\[
\psi^{M_2}\cirk\alpha_{\top,\ldots,\top}=\mj_\top\cirk\psi^{M_1}.
\]
This equation is the nullary case of $(\psi\alpha)$. Then the
supposition that $\delta^\str$ and $\sigma^\str$ are upward
preserved by $\top$ gives (VII).

We will now derive some more equations like $(c^m34)$, $(c^m234)$
and $(c^m24)$, which we need to prove that \aA\ is indeed a
symmetric monoidal category. First we have
\begin{tabbing}
\mbox{\hspace{2.5em}}$(c^m12)$\mbox{\hspace{2em}}\=
$c^m_{\top,\top,A_3,A_4}=(\sigma^\rts_{A_3}\kon
\sigma^\rts_{A_4})\cirk \sigma^\str_{A_3\kon
A_4}\cirk(\sigma^\str_\top\kon\mj_{A_3\kon A_4})$
\end{tabbing}
from $(\psi\sigma)$ and (VII), and from that we obtain
\begin{tabbing}
\mbox{\hspace{2.5em}}\=$(c^m34)$\mbox{\hspace{2em}}\=
$c^m_{A_1,A_2,\top,\top}$\kill \>$(c^m123)$\>
$c^m_{\top,\top,\top,A}=\mj_{(\top\kon\top)\kon(\top\kon A)}$
\end{tabbing}
(see how we obtained $(c^m234)$ above). From $(c^m12)$ and
$(c^mc^m)$ we derive
\begin{tabbing}
\mbox{\hspace{2.5em}}$(c^m13)$\mbox{\hspace{2em}}\=
$c^m_{\top,A_2,\top,A_4}=(\sigma^\rts_\top\kon \mj_{A_2\kon
A_4})\cirk \sigma^\rts_{A_2\kon
A_4}\cirk(\sigma^\str_{A_2}\kon\sigma^\str_{A_4})$.
\end{tabbing}

From $(\psi c^m)$, by substituting $\top$ for $A_1'$, $A_2$,
$A_2'$, $A_3$, $A_3'$ and $A_4$, together with $(c^m234)$,
$(c^m123)$ and $(c^mc^m)$ we easily obtain
\[
c^m_{A_1\kon\top,\top\kon\top,\top\kon\top,\top\kon A_4'}=
\mj_{((A_1\kon\top)\kon(\top\kon\top))\kon((\top\kon\top)\kon(\top\kon
A_4'))}.
\]
By precomposing and postcomposing with
\[
(\delta^\rts_{A_1}\kon\delta^\rts_\top)\kon(\delta^\rts_\top
\kon\sigma^\rts_{A_4'})
\]
and its inverse respectively, together with the naturality of
$c^m$ we obtain, after replacing $A_4'$ by $A_4$,
\begin{tabbing}
\mbox{\hspace{2.5em}}\=$(c^m13)$\mbox{\hspace{2em}}\=
$c^m_{\top,A_2,\top,A_4}=$\kill \>$(c^m23)$\mbox{\hspace{2em}}\>
$c^m_{A_1,\top,\top,A_4}=\mj_{(A_1\kon\top)\kon(\top\kon A_4)}$.
\end{tabbing}

We define in \aA\ the associativity arrows in the following
manner:
\begin{tabbing}
\mbox{\hspace{2.5em}}\=$b^\str_{A,B,C}\:$\=$=_{df} (\mj_{A\kon
B}\kon\sigma^\str_C)\cirk
c^m_{A,\top,B,C}\cirk(\delta^\rts_A\kon\mj_{B\kon C})\!:$\\[.5ex]
\`$A\kon(B\kon C)\str(A\kon B)\kon C$,\\[1ex]
\>$b^\rts_{A,B,C}$\>$=_{df} (\delta^\str_A\kon \mj_{B\kon C})\cirk
c^m_{A,B,\top,C}\cirk(\mj_{A\kon B}\kon\sigma^\rts_C)\!:$\\[.5ex]
\`$(A\kon B)\kon C\str A\kon(B\kon C)$.
\end{tabbing}
It is clear that $b^\str$ and $b^\rts$ are natural isomorphisms,
inverse to each other. Our purpose now is to derive in \aA\ Mac
Lane's pentagonal equation
\begin{tabbing}
\mbox{\hspace{2.5em}}\=$(b5)$\mbox{\hspace{2em}}\=$b^\str_{A\kon
B,C,D}\cirk b^\str_{A,B,C\kon D}=(b^\str_{A,B,C}\kon\mj_D)\cirk
b^\str_{A,B\kon C,D}\cirk(\mj_A\kon b^\str_{B,C,D})$.
\end{tabbing}

\prop{Proposition 1}{The equation $(b5)$ holds in \aA.}

\dkz With $LHS$ being the left-hand side of $(b5)$, we have
\begin{tabbing}
\mbox{\hspace{2.5em}}$LHS=(\mj_{(A\kon B)\kon
C}\kon\sigma^\str_D)\cirk c^m_{A\kon
B,\top,C,D}\cirk(\delta^\rts_{A\kon B}\kon\sigma^\str_{C\kon
D})\cirk c^m_{A,\top,B,C\kon
D}\cirk$\\[.5ex]
\`$\cirk (\delta^\rts_A\kon\mj_{B\kon (C\kon D)})$.
\end{tabbing}
By replacing $\delta^\rts_{A\kon B}$ and $\sigma^\str_{C\kon D}$
with the help of $(\psi\delta)$ and $(\psi\sigma)$, together with
bifunctorial and naturality equations, for
\begin{tabbing}
\mbox{\hspace{8em}}\=$f=((\mj_{A\kon
B}\kon\sigma^\str_C)\kon(\sigma^\str_D\cirk
(\delta^\str_\top\kon\sigma^\str_D))$,\\[1ex]
\>$g=((\delta^\rts_A\kon\delta^\rts_\top)\cirk\delta^\rts_A)\kon
(\delta^\rts_B\kon\mj_{C\kon D})$,
\end{tabbing}
we obtain
\begin{tabbing}
$LHS$\=$=f\cirk c^m_{A\kon B,\top\kon\top,\top\kon C,\top\kon
D}\cirk(c^m_{A,\top,B,\top}\kon c^m_{\top,\top,C,D})\cirk
c^m_{A\kon\top,\top\kon\top,B\kon\top,C\kon D}\cirk g$.
\end{tabbing}

With $RHS$ being the right-hand side of $(b5)$, together with
bifunctorial and naturality equations, and the isomorphism of
$\delta^\str_\top$, we obtain
\begin{tabbing}
\mbox{\hspace{2.5em}}$RHS=f\cirk
(c^m_{A,\top,B,C}\kon\mj_{(\top\kon\top)\kon(\top\kon D)})\cirk
c^m_{A\kon\top,\top\kon\top,B\kon C,\top\kon D}\cirk$\\[.5ex]
\`$\cirk(\mj_ {(A\kon\top)\kon(\top\kon\top)}\kon
c^m_{B,\top,C,D})\cirk g$.
\end{tabbing}
Then we use $(c^m123)$, $(c^m234)$ and $(\psi c^m)$ to obtain that
$LHS=RHS$.\qed

From the definition of $b^\str_{A,\top,C}$ with $(c^m23)$ we
obtain easily
\begin{tabbing}
\mbox{\hspace{2.5em}}$(b\delta\sigma)$\mbox{\hspace{2em}}
$b^\str_{A,\top,C}=\delta^\rts_A\kon\sigma^\str_C$.
\end{tabbing}
With that we have obtained in \aA\ all the equations of monoidal
categories; i.e., \aA\ is a monoidal category. Hence we can apply
Mac Lane's monoidal coherence theorem, which says that if
$f,g\!:A\str B$ are monoidal arrows, i.e.\ arrows definable in
terms of $\mj$, $b^\str$, $b^\rts$, $\delta^\str$, $\delta^\rts$,
$\sigma^\str$, $\sigma^\rts$, $\kon$ and $\cirk$, then ${f=g}$ in
monoidal categories (see \cite{ML63}, \cite{ML71}, Section VII.2,
or \cite{DP04}, Section 4.6).

We define in \aA\ the commutativity arrows in the following
manner:
\begin{tabbing}
\mbox{\hspace{2.5em}}\=$c_{A,B}=_{df}
(\sigma^\str_B\kon\delta^\str_A)\cirk
c^m_{\top,A,B,\top}\cirk(\sigma^\rts_A\kon\delta^\rts_B)\!:A\kon
B\str B\kon A$.
\end{tabbing}
It is clear that $c$ is a natural isomorphism, inverse to itself,
due to $(c^mc^m)$. Our purpose now is to derive in \aA\ Mac Lane's
hexagonal equation
\begin{tabbing}
\mbox{\hspace{2.5em}}\=$(bc)$\mbox{\hspace{2em}}\=$b^\rts_{C,A,B}\cirk
(c_{A,C}\kon\mj_B)\cirk b^\str_{A,C,B}\cirk(\mj_A\kon
c_{B,C})\cirk b^\rts_{A,B,C}=c_{A\kon B,C}$.
\end{tabbing}

\prop{Proposition 2}{The equation $(bc)$ holds in \aA.}

\dkz With $LHS$ being the left-hand side of $(bc)$ we have,
together with bifunctorial, naturality and isomorphism equations:
\begin{tabbing}
$LHS$\=$=f_1\cirk(c^m_{\top,A,C,\top}\kon\mj_{(\top\kon\top)\kon(B\kon\top)})
\cirk f_2\cirk(\mj_{((\top\kon A)\kon\top)\kon(\top\kon\top)}\kon
c^m_{\top,B,C,\top})\cirk f_3$
\end{tabbing}
for $f_1$, $f_2$ and $f_3$ monoidal arrows. Then by using
$(c^m23)$ and $(\psi c^m)$ we obtain
\[f_2\cirk(\mj_{((\top\kon
A)\kon\top)\kon(\top\kon\top)}\kon c^m_{\top,B,C,\top})\cirk
f_3=f'_2\cirk c^m_{(\top\kon A)\kon\top,\top\kon B,\top\kon
C,\top\kon\top}\cirk f'_3
\]
for $f'_2$ and $f'_3$ monoidal arrows. By previously derived
equations, together with bifunctorial and naturality equations, we
obtain
\begin{tabbing}
\mbox{\hspace{2.5em}}$LHS$\=$=f_1\cirk(c^m_{\top,A,C,\top}\kon
c^m_{\top,B,\top,\top}) \cirk c^m_{\top\kon A,\top\kon
B,C\kon\top,\top\kon\top}\cirk f''_3$
\end{tabbing}
for $f_1$ and $f''_3$ monoidal arrows. Then by using once more
$(\psi c^m)$, together with previously derived equations, we
obtain
\begin{tabbing}
\mbox{\hspace{2.5em}}$LHS$\=$=f'_1\cirk c^m_{\top\kon\top,A\kon
B,C\kon\top,\top\kon\top}\cirk f'''_3$
\end{tabbing}
for $f'_1$ and $f'''_3$ monoidal arrows. With the help of monoidal
coherence, together with bifunctorial, naturality and isomorphism
equations, we obtain $LHS=c_{A\kon B,C}$.\qed

With that we have obtained in \aA\ all the equations of symmetric
monoidal categories; i.e., \aA\ is a symmetric monoidal category.
To ascertain that we have obtained not more than the equations of
symmetric monoidal categories it remains to derive the definition
of $c^m$ in terms of $b^\str$, $b^\rts$ and $c$---namely, the
equation
\begin{tabbing}
\mbox{\hspace{2.5em}}\=$(c^mbc)$\mbox{\hspace{2em}}\=$c^m_{A,B,C,D}=
b^\str_{A,C,B\kon
D}\cirk(\mj_A\kon(b^\rts_{C,B,D}\cirk(c_{B,C}\kon\mj_D)\cirk
b^\str_{B,C,D}))\cirk$\\[1ex]
\`$\cirk b^\rts_{A,B,C\kon D}$.
\end{tabbing}
\vspace{-2ex} \prop{Proposition 3}{The equation $(c^mbc)$ holds in
\aA.}

\dkz With bifunctorial, naturality and isomorphism equations we
have
\begin{tabbing}
\mbox{\hspace{2.5em}}$(c^mbc)$\mbox{\hspace{2em}}\=$c^m_{A,B,C,D}$\kill
\>$c^m_{A,B,C,D}=f_1\cirk c^m_{A\kon\top,\top\kon B,\top\kon
C,\top\kon D}\cirk f_2$
\end{tabbing}
for $f_1$ and $f_2$ monoidal arrows. By using $(\psi c^m)$ and
previously derived equations we obtain
\begin{tabbing}
\mbox{\hspace{2.5em}}$(c^mbc)$\mbox{\hspace{2em}}\=$c^m_{A,B,C,D}$\kill
\>$c^m_{A,B,C,D}=f'_1\cirk(\mj_{(A\kon\top)\kon(\top\kon\top)}\kon
c^m_{\top,B,C,D})\cirk f'_2$
\end{tabbing}
for $f'_1$ and $f'_2$ monoidal arrows. With bifunctorial,
naturality and isomorphism equations we have
\begin{tabbing}
\mbox{\hspace{2.5em}}$(c^mbc)$\mbox{\hspace{2em}}\=$c^m_{A,B,C,D}$\kill
\>$c^m_{\top,B,C,\top}=f_3\cirk c^m_{\top\kon\top,B
\kon\top,C\kon\top,\top\kon D}\cirk f_4$
\end{tabbing}
for $f_3$ and $f_4$ monoidal arrows. By using once more $(\psi
c^m)$, together with previously derived equations, we obtain
\begin{tabbing}
\mbox{\hspace{2.5em}}\=$c^m_{A,B,C,D}$\kill
\>$c^m_{A,B,C,D}=f''_1\cirk(\mj_{(A\kon\top)\kon(\top\kon\top)}\kon
(c^m_{\top,B,C,\top}\kon\mj_{(\top\kon\top)\kon(\top\kon
D)}))\cirk f''_2$
\end{tabbing}
for $f''_1$ and $f''_2$ monoidal arrows. From this equation with
the help of monoidal coherence, together with bifunctorial,
naturality and isomorphism equations, we obtain $(c^mbc)$.\qed

With $c^m$, instead of $b^\str$, $b^\rts$ and $c$, added to $\mj$,
$\delta^\str$, $\delta^\rts$, $\sigma^\str$, $\sigma^\rts$, $\kon$
and $\cirk$ we have an alternative language for symmetric monoidal
categories. In this alternative language symmetric monoidal
categories are defined by the categorial, bifunctorial, naturality
and isomorphism equations to which we add the equations $(\psi
c^m)$, $(\psi\delta)$, $(\psi\sigma)$ and (VII). All these
equations hold in symmetric monoidal categories defined in a
standard way when $c^m$ is defined in terms of $b^\str$, $b^\rts$
and $c$ according to $(c^mbc)$. In the alternative definition of
symmetric monoidal category, to define $b^\str$, $b^\rts$ and $c$
it is essential to have besides $c^m$ the unit object $\top$ and
the associated arrows.

\section{\large\bf The equation $(\psi c^m)$ and the Yang-Baxter equation}
As we announced in Section~1, the present section is a brief
comment upon the equation $(\psi c^m)$. We will not go into
details, but just point to references where the subject is treated
more exhaustively.

The key equation that delivers Mac Lane's pentagonal and hexagonal
equations of symmetric monoidal categories is the equation $(\psi
c^m)$. This equation is closely related to the Yang-Baxter
equation, forms of which are found in symmetric and braid groups
(see the equation (YB) in Section~5). In symmetric and braided
monoidal categories where the $b^\str$ and $b^\rts$ arrows are
identity arrows this equation becomes
\[
(c_{B,C}\kon\mj_A)\cirk(\mj_B\kon c_{A,C})\cirk(c_{A,B}\kon\mj_C)=
(\mj_C\kon c_{A,B})\cirk(c_{A,C}\kon\mj_B)\cirk(\mj_A\kon
c_{B,C}).
\]

Let us draw graphs like those one finds in Mac Lane's coherence
result for symmetric monoidal categories (see \cite{ML63}, or
\cite{DP04}, Section 2.9 and Chapter~5; cf.\ the functor $G$ in
Section~5). The graphs that correspond to the two sides of the
Yang-Baxter equation above are
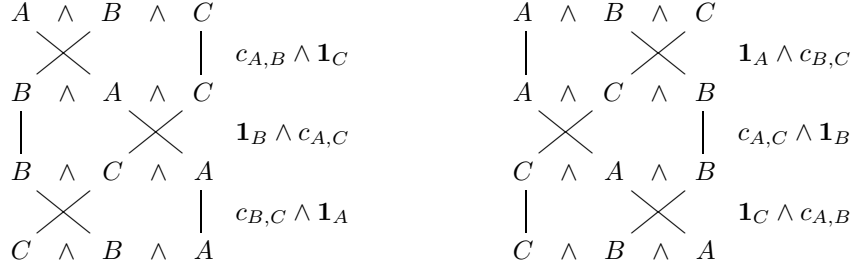
\begin{center}
\begin{picture}(340,110)
\put(10,10){$C$\quad$\kon$\quad$B$\quad$\kon$\quad$A$}

\put(10,40){$B$\quad$\kon$\quad$C$\quad$\kon$\quad$A$}

\put(10,70){$B$\quad$\kon$\quad$A$\quad$\kon$\quad$C$}

\put(10,100){$A$\quad$\kon$\quad$B$\quad$\kon$\quad$C$}

\put(200,10){$C$\quad$\kon$\quad$B$\quad$\kon$\quad$A$}

\put(200,40){$C$\quad$\kon$\quad$A$\quad$\kon$\quad$B$}

\put(200,70){$A$\quad$\kon$\quad$C$\quad$\kon$\quad$B$}

\put(200,100){$A$\quad$\kon$\quad$B$\quad$\kon$\quad$C$}

\put(95,25){$c_{B,C}\kon\mj_A$}

\put(95,55){$\mj_B\kon c_{A,C}$}

\put(95,85){$c_{A,B}\kon\mj_C$}

\put(285,25){$\mj_C\kon c_{A,B}$}

\put(285,55){$c_{A,C}\kon\mj_B$}

\put(285,85){$\mj_A\kon c_{B,C}$}

\put(20,20){\line(5,4){20}}

\put(20,80){\line(5,4){20}}

\put(55,50){\line(5,4){20}}

\put(245,20){\line(5,4){20}}

\put(245,80){\line(5,4){20}}

\put(210,50){\line(5,4){20}}

\put(40,20){\line(-5,4){20}}

\put(40,80){\line(-5,4){20}}

\put(75,50){\line(-5,4){20}}

\put(265,20){\line(-5,4){20}}

\put(265,80){\line(-5,4){20}}

\put(230,50){\line(-5,4){20}}

\put(82,20){\line(0,1){16}}

\put(82,80){\line(0,1){16}}

\put(14,50){\line(0,1){16}}

\put(205,20){\line(0,1){16}}

\put(205,80){\line(0,1){16}}

\put(272,50){\line(0,1){16}}
\end{picture}
\end{center}
\vspace{-2ex} Passing from one of these graphs to the other is
analogous to the third Reidemeister move of knot theory.

The scheme of the graph that corresponds to $c^m_{A,B,C,D}$ is
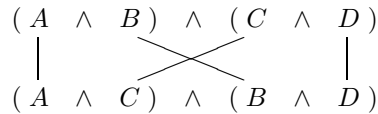
\begin{center}
\begin{picture}(170,50)
\put(10,10){$(\;A$\quad$\kon$\quad$C\;)$\quad$\kon$\quad$(\;B$\quad$\kon$\quad$D\;)$}

\put(10,40){$(\;A$\quad$\kon$\quad$B\;)$\quad$\kon$\quad$(\;C$\quad$\kon$\quad$D\;)$}

\put(58,20){\line(5,2){40}}

\put(98,20){\line(-5,2){40}}

\put(20,20){\line(0,1){16}}

\put(137,20){\line(0,1){16}}
\end{picture}
\end{center}
\vspace{-2ex} (as dictated by the definition $(c^mbc)$ of $c^m$ in
terms of $b^\str$, $b^\rts$ and $c$; see the functor $G$ in
Section~5), and here are the graphs that correspond to the two
sides of $(\psi c^m)$, with some lines drawn dotted for reasons
given below,
\begin{center}
\begin{picture}(330,108)
\put(0,100){$((A_1\kon A_1')\kon(A_2\kon A_2'))\kon((A_3\kon
A_3')\kon(A_4\kon A_4'))$}

\put(0,70){$((A_1\kon A_1')\kon(A_3\kon A_3'))\kon((A_2\kon
A_2')\kon(A_4\kon A_4'))$}

\put(225,88){\small $c^m_{A_1\kon A_1',A_2\kon A_2',A_3\kon
A_3',A_4\kon A_4'}$}

\put(75,80){\line(3,1){45}}

\multiput(98,80)(3,1){16}{\circle*{.3}}

\put(120,80){\line(-3,1){45}}

\multiput(143,80)(-3,1){16}{\circle*{.3}}

\multiput(12,80)(0,6){3}{\line(0,1){4}}

\put(35,80){\line(0,1){16}}

\multiput(182,80)(0,3){6}{\circle*{.3}}

\multiput(205,80)(0,6){3}{\line(0,1){4}}


\put(38,50){\line(3,2){23}}

\put(63,50){\line(-3,2){23}}

\multiput(154,50)(3,2){9}{\circle*{.3}}

\multiput(179,50)(-3,2){9}{\circle*{.3}}

\multiput(12,50)(0,6){3}{\line(0,1){4}}

\multiput(91,50)(0,3){6}{\circle*{.3}}

\multiput(205,50)(0,6){3}{\line(0,1){4}}

\put(127,50){\line(0,1){16}}

\put(0,40){$((A_1\kon A_3)\kon(A_1'\kon A_3'))\kon((A_2\kon
A_4)\kon(A_2'\kon A_4'))$}

\put(225,58){\small $c^m_{A_1,A_1',A_3,A_3'}\kon
c^m_{A_2,A_2',A_4,A_4'}$}


\put(75,20){\line(3,1){45}}

\multiput(98,20)(3,1){16}{\circle*{.3}}

\put(120,20){\line(-3,1){45}}

\multiput(143,20)(-3,1){16}{\circle*{.3}}

\multiput(12,20)(0,6){3}{\line(0,1){4}}

\put(35,20){\line(0,1){16}}

\multiput(182,20)(0,3){6}{\circle*{.3}}

\multiput(205,20)(0,6){3}{\line(0,1){4}}

\put(0,10){$((A_1\kon A_3)\kon(A_2\kon A_4))\kon((A_1'\kon
A_3')\kon(A_2'\kon A_4'))$}

\put(225,28){\small $c^m_{A_1\kon A_3,A_1'\kon A_3',A_2\kon
A_4,A_2'\kon A_4'}$}
\end{picture}
\end{center}
\begin{center}
\begin{picture}(330,108)
\put(0,100){$((A_1\kon A_1')\kon(A_2\kon A_2'))\kon((A_3\kon
A_3')\kon(A_4\kon A_4'))$}

\put(0,70){$((A_1\kon A_2)\kon(A_1'\kon A_2'))\kon((A_3\kon
A_4)\kon(A_3'\kon A_4'))$}

\put(225,88){\small $c^m_{A_1,A_1',A_2,A_2'}\kon
c^m_{A_3,A_3',A_4,A_4'}$}

\put(38,80){\line(3,2){23}}

\put(63,80){\line(-3,2){23}}

\multiput(154,80)(3,2){9}{\circle*{.3}}

\multiput(179,80)(-3,2){9}{\circle*{.3}}

\multiput(12,80)(0,6){3}{\line(0,1){4}}

\multiput(91,80)(0,3){6}{\circle*{.3}}

\multiput(205,80)(0,6){3}{\line(0,1){4}}

\put(127,80){\line(0,1){16}}


\put(0,40){$((A_1\kon A_2)\kon(A_3\kon A_4))\kon((A_1'\kon
A_2')\kon(A_3'\kon A_4'))$}

\put(225,58){\small $c^m_{A_1\kon A_2,A_1'\kon A_2',A_3\kon
A_4,A_3'\kon A_4'}$}

\put(75,50){\line(3,1){45}}

\multiput(98,50)(3,1){16}{\circle*{.3}}

\put(120,50){\line(-3,1){45}}

\multiput(143,50)(-3,1){16}{\circle*{.3}}

\multiput(12,50)(0,6){3}{\line(0,1){4}}

\put(35,50){\line(0,1){16}}

\multiput(182,50)(0,3){6}{\circle*{.3}}

\multiput(205,50)(0,6){3}{\line(0,1){4}}


\put(38,20){\line(3,2){23}}

\put(63,20){\line(-3,2){23}}

\multiput(154,20)(3,2){9}{\circle*{.3}}

\multiput(179,20)(-3,2){9}{\circle*{.3}}

\multiput(12,20)(0,6){3}{\line(0,1){4}}

\put(127,20){\line(0,1){16}}

\multiput(205,20)(0,6){3}{\line(0,1){4}}

\multiput(91,20)(0,3){6}{\circle*{.3}}

\put(0,10){$((A_1\kon A_3)\kon(A_2\kon A_4))\kon((A_1'\kon
A_3')\kon(A_2'\kon A_4'))$}

\put(225,28){\small $c^m_{A_1,A_2,A_3,A_4}\kon
c^m_{A_1',A_2',A_3',A_4'}$}
\end{picture}
\end{center}
\vspace{-2ex}

Consider in each of these two graphs the subgraph $G_1$ made of
solid lines, involving $A_1'$, $A_2$ and $A_3$. This subgraph is
interlaced with the subgraph $G_2$ made of finely dotted lines,
involving $A_2'$, $A_3'$ and $A_4$. (The principle of this
interlacing, and the connection of $G_1$ and $G_2$ with $G^3_1$
and $G^3_2$ respectively, is explained in Section~5.) The
subgraphs $G_1$ and $G_2$ of each side of $(\psi c^m)$ are like
the graphs of the two sides of the Yang-Baxter equation. On the
other hand, $G_1$ of one side of $(\psi c^m)$ and $G_1$ of the
other side are also like the graphs of the two sides of the
Yang-Baxter equation, and analogously for $G_2$.

In the Brauerian representation (see \cite{DP03a}, \cite{DP03b},
\cite{DP05} and references therein; this representation stems from
\cite{B37}), the graph $B$(X), which corresponds to medial
commutativity:
\begin{center}
\begin{picture}(130,48)
\put(0,40){$\;\;\;$00$\;\;\;\;\;\;$01$\;\;\;\;\;$10$\;\;\;\;\;\;$11}

\put(0,10){$\;\;\;$00$\;\;\;\;\;\;$01$\;\;\;\;\;$10$\;\;\;\;\;\;$11}

\put(38,20){\line(3,2){23}}

\put(63,20){\line(-3,2){23}}

\put(12,20){\line(0,1){16}}

\put(91,20){\line(0,1){16}}
\end{picture}
\end{center}
\vspace{-2ex} stands for the graph X:
\begin{center}
\begin{picture}(50,20)

\put(0,5){\line(5,4){20}}

\put(20,5){\line(-5,4){20}}
\end{picture}
\end{center}
\vspace{-2ex} The binary numbering of places in $B$(X) is required
by the representation. (We can read 0 as ``left'' and 1 as
``right''.) Let us say that a graph can be \emph{fitted} between
two 0-1 sequences, one above the other, when the lines of the
graph connect two occurrences of 0 or two occurrences of 1, and no
occurrence of 0 with an occurrence of 1. The graph X can be fitted
between 00 and 00, between 01 and 10, between 10 and 01, and
between 11 and 11. In $B$(X) lines join places between which X can
be fitted. This is the gist of the Brauerian representation. Note
that this representation is an isomorphic representation.

In the Brauerian representation, the graphs of the two sides of
the Yang-Baxter equation are represented by graphs that are
exactly like the graphs of the two sides of $(\psi c^m)$. To
obtain the latter graphs we compose copies of the graph:
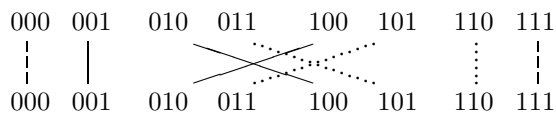
\begin{center}
\begin{picture}(230,48)
\put(0,40){$\;\;$000$\;\;\;$001$\;\;\;\;\;$010$\;\;\;\;$011$\;\;\;\;\;\;\;$100$\;\;\;\;$101$\;\;\;\;\;$110$\;\;\;$111}

\put(0,10){$\;\;$000$\;\;\;$001$\;\;\;\;\;$010$\;\;\;\;$011$\;\;\;\;\;\;\;$100$\;\;\;\;$101$\;\;\;\;\;$110$\;\;\;$111}

\put(75,20){\line(3,1){45}}

\multiput(98,20)(3,1){16}{\circle*{.3}}

\put(120,20){\line(-3,1){45}}

\multiput(143,20)(-3,1){16}{\circle*{.3}}

\multiput(12,20)(0,6){3}{\line(0,1){4}}

\put(35,20){\line(0,1){16}}

\multiput(182,20)(0,3){6}{\circle*{.3}}

\multiput(205,20)(0,6){3}{\line(0,1){4}}
\end{picture}
\end{center}
\vspace{-2ex}which in the Brauerian representation stands for the
graph:
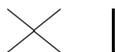
\begin{center}
\begin{picture}(50,20)

\put(0,5){\line(5,4){20}}

\put(20,5){\line(-5,4){20}}

\put(40,5){\line(0,1){16}}
\end{picture}
\end{center}
\vspace{-2ex} with copies of the graph:
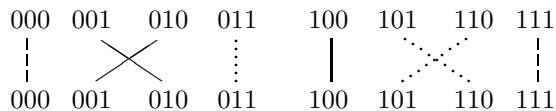
\begin{center}
\begin{picture}(230,48)
\put(0,40){$\;\;$000$\;\;\;$001$\;\;\;\;\;$010$\;\;\;\;$011$\;\;\;\;\;\;\;$100$\;\;\;\;$101$\;\;\;\;\;$110$\;\;\;$111}

\put(0,10){$\;\;$000$\;\;\;$001$\;\;\;\;\;$010$\;\;\;\;$011$\;\;\;\;\;\;\;$100$\;\;\;\;$101$\;\;\;\;\;$110$\;\;\;$111}

\put(38,20){\line(3,2){23}}

\put(63,20){\line(-3,2){23}}

\multiput(154,20)(3,2){9}{\circle*{.3}}

\multiput(179,20)(-3,2){9}{\circle*{.3}}

\multiput(12,20)(0,6){3}{\line(0,1){4}}

\put(127,20){\line(0,1){16}}

\multiput(205,20)(0,6){3}{\line(0,1){4}}

\multiput(91,20)(0,3){6}{\circle*{.3}}
\end{picture}
\end{center}
\vspace{-2ex} which in the Brauerian representation stands for the
graph:
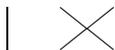
\begin{center}
\begin{picture}(50,20)

\put(20,5){\line(5,4){20}}

\put(40,5){\line(-5,4){20}}

\put(0,5){\line(0,1){16}}
\end{picture}
\end{center}
Note that, with the binary numbering of places, the subgraphs
$G_1$ involve places with just one occurrence of 1, and the
subgraphs $G_2$ places with just two occurrences of 1. The
remaining straight-line subgraphs at the far left and the far
right involve places with 0 and 3 occurrences of 1 respectively.
(This is related to the subject of Section~5, and the subject of
\cite{JK83}, Chapter 2.)

In \cite{DP06} it was shown how the pentagon may be conceived as a
degenerate Yang-Baxter hexagon, in which two vertices are joined
into one and one side collapses. This collapse was previously
described geometrically and combinatorially in \cite{T97}.

In deriving the pentagonal equation $(b5)$ from $(\psi c^m)$ in
the proof of Proposition~1 of the preceding section we also
obtained a collapse of one side of the hexagon of $(\psi c^m)$,
which as we said above corresponds through the Brauerian
representation to the Yang-Baxter hexagon. The side collapsed in
the derivation of $(b5)$ and the side collapsed in \cite{DP06} are
the same.

Another connection of $(\psi c^m)$ with the Yang-Baxter equation
is made manifest by considering the permutations of the indices
$i$, $j$ and $k$ in the hexagonal interchange equation of
\cite{BFSV03} (end of Definition 1.7). The equation $(\psi c^m)$
is obtained from this hexagonal interchange equation by deleting
these three indices (as we mentioned in the preceding section).

\section{\large\bf Associative categories}
Associative categories are monoidal categories without unit
object. We have in them, namely, a biendofunctor $\kon$ and the
natural isomorphism $b^\str$ for which Mac Lane's pentagonal
equation $(b5)$ holds. Our purpose now is to explain this notion
in the same style as we explained the notion of symmetric monoidal
category in Section~2. At the end of the section we consider
briefly monoidal categories, i.e.\ associative categories with
unit object, in the same light.

To make the connection with Section~2 precise, we will first
introduce a generalization that covers both the previous matter
and what we need for associative categories. This will complicate
notation considerably---there will be many iterated indices---but,
nevertheless, the subject is not very difficult.

In this section we talk only of shapes without $\top$, whose arity
is hence at least $1$. If $N$ is a shape of arity $n$ and $M$ a
shape of arity $m$, let $N_M^i$ be the shape obtained by replacing
the $i$-th occurrence of \koc\ in $N$ (counting from the left) by
the shape $M$. We have, of course, $1\leq i\leq n$, and the arity
of $N_M^i$ is $n\pl m\mn 1$. The shape $N_\Box^i$ is $N$ and
$\koc_M^1$ is $M$. For $n_1$ and $n_2$ arities of the shapes $N_1$
and $N_2$ respectively, we have
\[
(N_1\kon N_2)_M^i=\left\{
\begin{array}{ll}
(N_1)_M^i\kon N_2 & \mbox{\rm for }1\leq i\leq n_1
\\[1.5ex]
N_1\kon(N_2)_M^{i-n_1} & \mbox{\rm for }n_1\pl 1\leq i\leq n_1\pl
n_2.
\end{array}
\right.
\]

The \emph{split} $[N_M^i]^S$ \emph{of} $N_M^i$ \emph{through} $M$
is a pair of shapes $([N_M^i]^L,[N_M^i]^R)$ defined inductively as
follows:
\[
[\koc_M^1]^S=([\koc_M^1]^L,[\koc_M^1]^R)=(M,M),
\]
\[
[(N_1\kon N_2)_M^i]^S=\left\{
\begin{array}{ll}
([(N_1)_M^i]^L,[(N_1)_M^i]^R\kon N_2) & \mbox{\rm for }1\leq i\leq
n_1
\\[1.5ex]
(N_1\kon [(N_2)_M^{i-n_1}]^L,[(N_2)_M^{i-n_1}]^R) & \mbox{\rm for
}n_1\pl 1\leq i\leq n_1\pl n_2.
\end{array}
\right.
\]
The shape $M$ is nested on the extreme right of $[N_M^i]^L$ and on
the extreme left of $[N_M^i]^R$.

For $F$ a biendofunctor of \aA, we will define below the natural
transformation $\psi^{N,i,M}$, whose components are the arrows of
\aA
\begin{tabbing}
\mbox{\hspace{2.5em}}\=$\psi^{N,i,M}_{\vec{B}_{i\mn 1},
\scriptsize\overrightarrow{A_m, A_m'},\vec{C}_{n\mn i}}\!\!\!:
N_M^i(\vec{B}_{i-1},\overrightarrow{F(A_m,A_m')},\vec{C}_{n-i})\str$
\\*[.5ex]
\`$F([N_M^i]^L(\vec{B}_{i-1},\vec{A}_m),[N_M^i]^R(\vec{A}_m',\vec{C}_{n-i}))$.
\end{tabbing}
Note that if $i=1$, then there are no $B$'s, and if $i=n$, then
there are no $C$'s. Note also that the formula
\[N_M^i(\vec{B}_{i-1},\overrightarrow{F(A_m,A_m')},\vec{C}_{n-i})\quad
\mbox{\rm{ is }}\quad
N(\vec{B}_{i-1},M(\overrightarrow{F(A_m,A_m')}),\vec{C}_{n-i}).\]

By the definition of split we have
\[
\begin{array}{l}
[(\koc\kon\koc)_\Box^1]^S=([\koc_\Box^1]^L,[\koc_\Box^1]^R\kon\koc)=
(\koc,\koc\kon\koc),
\\[1.5ex]
[(\koc\kon\koc)_\Box^2]^S=(\koc\kon[\koc_\Box^1]^L,[\koc_\Box^1]^R)=
(\koc\kon\koc,\koc),
\end{array}
\]
so that for $F$ being $\kon$ we have
\[
\begin{array}{l}
\psi_{A_1,A_1',C}^{\Box\kon\Box,1,\Box}\!:(A_1\kon A_1')\kon C\str
A_1\kon(A_1'\kon C),
\\[1.5ex]
\psi_{B,A_1,A_1'}^{\Box\kon\Box,2,\Box}\!:B\kon(A_1\kon
A_1')\str(B\kon A_1)\kon A_1'.
\end{array}
\]
We write $b_{A,B,C}^\rts$ for $\psi_{A,B,C}^{\Box\kon\Box,1,\Box}$
and $b_{A,B,C}^\str$ for $\psi_{A,B,C}^{\Box\kon\Box,2,\Box}$.

Now comes the inductive definition of $\psi^{N,i,M}$ in terms of
$\psi^M$ of Section~2 and the $b^\rts$ and $b^\str$ arrows:

\vspace{1ex}

\noindent\hspace{2em}for $i=n=1$,
\[
\psi_{\scriptsize\overrightarrow{A_m,A_m'}}^{\Box,1,M}=
\psi^M_{\scriptsize\overrightarrow{A_m,A_m'}};
\]
\hspace{2em}for $1\leq i\leq n_1$,

\begin{tabbing}
\mbox{\hspace{0em}}\=$\psi_{\vec{B}_{i-1},\scriptsize\overrightarrow{A_m,A_m'}
\vec{C}_{n_1-i},\vec{D}_{n_2}}^{N_1\kon
N_2,i,M}=b^\rts_{[(N_1)_M^i]^L(\vec{B}_{i-1},\vec{A}_m)\:,\:
[(N_1)_M^i]^R(\vec{A}'_m,
\vec{C}_{n_1-i})\:,\:N_2(\vec{D}_{n_2})}$
\\*[.5ex]
\`$\cirk(\psi_{\vec{B}_{i-1},\scriptsize\overrightarrow{A_m,A_m'},
\vec{C}_{n_1-i}}^{N_1,i,M} \kon\mj_{N_2(\vec{D}_{n_2})})$;
\end{tabbing}
\hspace{2em}for $n_1\pl 1\leq i\leq n_1\pl n_2$,

\begin{tabbing}
\mbox{\hspace{0em}}\=$\psi_{\vec{D}_{n_1},\vec{B}_{i-n_1-1},
\scriptsize\overrightarrow{A_m,A_m'} \vec{C}_{n_2-i}}^{N_1\kon
N_2,i,M}\!\!\!\!\!\!\!\!=b^\str_{N_1(\vec{D}_{n_1})\:,\:
[(N_2)_M^{i-n_1}]^L(\vec{B}_{i-n_1-1},\vec{A}_m)\:,\:
[(N_2)_M^{i-n_1}]^R(\vec{A}'_m,\vec{C}_{n_2-i})}$
\\*[.5ex]
\`$\cirk(\mj_{N_1(\vec{D}_{n_1})}\kon
\psi_{\vec{B}_{i-n_1-1},\scriptsize\overrightarrow{A_m,A_m'},\vec{C}_{n_2-i}}^{N_2,i-n_1,M}
)$.
\end{tabbing}
Note that there are no occurrences of $c^m$ in
$\psi_{\vec{B}_{i\mn 1},A,A',\vec{C}_{n\mn i}}^{N,i,\Box}$.

\vspace{1ex}

Let $N_1$ and $N_2$ be shapes of arity $n$, while $M_1$ and $M_2$
are shapes of arity $m$, and let $\pi$ be a permutation of
$\{1,\ldots,n\pl m\mn 1\}$ that keeps fixed the first $i\mn 1$
elements and the last $n\mn i$ elements. Let $\alpha$ be a natural
transformation from the functor defined by the shape
$(N_1)_{M_1}^i$ to the functor $((N_2)_{M_2}^i)^{\pi}$ defined by
the shape $(N_2)_{M_2}^i$ and the permutation $\pi$ (see
Section~2). Let $\pi^L$ be $\pi$ restricted to $\{1,\ldots,i\mn
1\pl m\}$ and let $\pi^R$ be $\pi$ restricted to $\{i\pl
1,\ldots,n\pl m\mn 1\}$.

Let $\alpha^L$ be a natural transformation from
$[(N_1)_{M_1}^i]^L$ to $([(N_2)_{M_2}^i]^L)^{\pi^L}$ and
$\alpha^R$ a natural transformation from $[(N_1)_{M_1}^i]^R$ to
$([(N_2)_{M_2}^i]^R)^{\pi^R}$. We say that $\alpha$ is
$i$-\emph{lifted to} $\alpha^L\kon\alpha^R$ when diagrams of the
following form commute in \aA:

\begin{center}
\begin{picture}(300,130)
\put(-21,50){\small
$[(N_1)^i_{M_1}]^L(\vec{B}_{i-1},\vec{A}_m)\kon
[(N_1)^i_{M_1}]^R(\vec{A}_m',\vec{C}_{n-i})$}

\put(80,5){\small
$([(N_2)^i_{M_2}]^L)^{\pi^L}(\vec{B}_{i-1},\vec{A}_{\pi(m)})\kon
([(N_2)^i_{M_2}]^R)^{\pi^R}(\vec{A}_{\pi(m)}',\vec{C}_{n-i})$}

\put(-5,120){\small
$(N_1)^i_{M_1}(\vec{B}_{i-1},\overrightarrow{A_m\kon
A_m'},\vec{C}_{n-i})$}

\put(175,75){\small
$(N_2)^i_{M_2}(\vec{B}_{i-1},\overrightarrow{A_{\pi(m)}\kon
A_{\pi(m)}'},\vec{C}_{n-i})$}

\put(143,109){\small $\alpha_{\vec{B}_{i\mn
1},\scriptsize\overrightarrow{A_m\kon A_m'},\vec{C}_{n\mn i}}$}

\put(110,110){\vector(4,-1){80}}

\put(57,28){\small $\alpha^L_{\vec{B}_{i\mn
1},\vec{A}_m}\kon\alpha^R_{\vec{A}_m', \vec{C}_{n\mn i}}$}

\put(110,40){\vector(4,-1){80}}

\put(12,90){\small $\psi^{N_1,i,M_1}_{\vec{B}_{i\mn 1},
\scriptsize\overrightarrow{A_m,A_m'},\vec{C}_{n\mn i}}$}

\put(90,110){\vector(0,-1){45}}

\put(210,68){\vector(0,-1){45}}

\put(215,48){\small $\psi^{N_2,i,M_2}_{\vec{B}_{i\mn 1},
\scriptsize\overrightarrow{A_{\pi(m)},A_{\pi(m)}'},\vec{C}_{n\mn
i}}$}
\end{picture}
\end{center}
i.e., we have in \aA\ the equation
\begin{tabbing}
\mbox{\hspace{2.5em}}$(\psi\alpha\ast)$\mbox{\hspace{2em}}
$\psi^{N_2,i,M_2}_{\vec{B}_{i\mn 1},
\scriptsize\overrightarrow{A_{\pi(m)},A_{\pi(m)}'},\vec{C}_{n\mn
i}}\!\!\cirk \alpha_{\vec{B}_{i\mn
1},\scriptsize\overrightarrow{A_m\kon
A_m'},\vec{C}_{n\mn i}}=$\\[1ex]
\`$(\alpha^L_{\vec{B}_{i\mn 1},\vec{A}_m}\kon\alpha^R_{\vec{A}_m',
\vec{C}_{n\mn i}})\cirk \psi^{N_1,i,M_1}_{\vec{B}_{i\mn 1},
\scriptsize\overrightarrow{A_m,A_m'},\vec{C}_{n\mn i}}$.
\end{tabbing}
The equation $(\psi\alpha)$ of Section~2 is a particular instance
of $(\psi\alpha*)$ when $N_1$ and $N_2$ are both \koc, $i=n=1$,
and $\alpha^L$ and $\alpha^R$ are both $\alpha$. So upward
preservation is a particular case of $i$-lifting.

We will make now a number of suppositions that will make out of
\aA\ at the end an associative category. The first supposition is
the same as in Section~2:
\begin{quote}
(I)\quad $\kon$ is a biendofunctor of the category \aA.
\end{quote}

We say that $\kon$ is $M$-\emph{lifting} in a shape $N$ of arity
$n$ when for every $i\in\{1,\ldots,n\}$ there is a natural
transformation $\psi^{N,i,M}$. In particular, $\kon$ is
\koc-lifting in a shape $N$ when for every $i\in\{1,\ldots,n\}$
there is a natural transformation whose components are the arrows
of \aA
\[
\psi_{\vec{B}_{i\mn 1},A,A',\vec{C}_{n\mn
i}}^{N,i,\Box}\!\!\!:N(\vec{B}_{i-1},A\kon A',\vec{C}_{n-
i})\str([N_\Box^i]^L(\vec{B}_{i-1},A)\kon[N_\Box^i]^R(A',
\vec{C}_{n-i})).
\]

The second supposition we make in order to make out of \aA\ an
associative category is the following:
\begin{quote}
(II)\quad $\kon$ is \koc-lifting in $\koc\kon\koc$.
\end{quote}
This implies that we have in \aA\ the natural transformations
$b^\str$ and $b^\rts$, and, as we saw above, this implies that
$\kon$ is \koc-lifting in every shape $N$ (which, as we assumed at
the beginning of the section, is without $\top$).

Our third supposition is the following:
\begin{quote}
(III)\quad $b^\str$ is 3-lifted to $b^\str\kon\mj$.
\end{quote}
Here $\alpha$ and $\alpha^L$ are both the natural transformation
$b^\str$ from $\koc\kon(\koc\kon\koc)$ to
$(\koc\kon\koc)\kon\koc$, while $\alpha^R$ is $\mj$ understood as
a natural transformation from \koc\ to \koc. Supposition (III)
means that in \aA\ the following instance of $(\psi\alpha*)$
commutes
\begin{center}
\begin{picture}(300,130)
\put(-23,50){\small
$[(\koc\kon(\koc\kon\koc))_\Box^3]^L(B_1,B_2,A)\kon
[(\koc\kon(\koc\kon\koc))_\Box^3]^R(A')$}

\put(100,5){\small
$[((\koc\kon\koc)\kon\koc)_\Box^3]^L(B_1,B_2,A)\kon
[((\koc\kon\koc)\kon\koc)_\Box^3]^R(A')$}

\put(-10,120){\small $(\koc\kon(\koc\kon\koc))_\Box^3
(B_1,B_2,A\kon A')$}

\put(180,75){\small $((\koc\kon\koc)\kon\koc)_\Box^3
(B_1,B_2,A\kon A')$}

\put(143,109){\small $b^\str_{B_1,B_2,A\kon A'}$}

\put(110,110){\vector(4,-1){80}}

\put(73,25){\small $b^\str_{B_1,B_2,A}\kon\mj_{A'}$}

\put(110,40){\vector(4,-1){80}}

\put(40,90){\small
$\psi^{\Box\kon(\Box\kon\Box),3,\Box}_{B_1,B_2,A,A'}$}

\put(97,110){\vector(0,-1){45}}

\put(220,68){\vector(0,-1){45}}

\put(225,48){\small
$\psi^{(\Box\kon\Box)\kon\Box,3,\Box}_{B_1,B_2,A,A'}$}
\end{picture}
\end{center}
which amounts to
\begin{center}
\begin{picture}(250,150)
\put(0,10){$(B_1\kon(B_2\kon A))\kon A'$}

\put(0,70){$B_1\kon((B_2\kon A)\kon A')$}

\put(0,130){$B_1\kon(B_2\kon(A\kon A'))$}

\put(180,10){$((B_1\kon B_2)\kon A)\kon A'$}

\put(180,130){$(B_1\kon B_2)\kon(A\kon A')$}

\put(105,4){$b^\str_{B_1,B_2,A}\kon\mj_{A'}$}

\put(100,15){\vector(1,0){70}}

\put(110,143){$b^\str_{B_1,B_2,A\kon A'}$}

\put(100,135){\vector(1,0){70}}

\put(4,45){$b^\str_{B_1,B_2\kon A,A'}$}

\put(55,60){\vector(0,-1){30}}

\put(-8,105){$\mj_{B_1}\kon b^\str_{B_2,A,A'}$}

\put(55,120){\vector(0,-1){30}}

\put(225,70){$b^\str_{B_1\kon B_2,A,A'}$}

\put(220,120){\vector(0,-1){90}}

\end{picture}
\end{center}
i.e.\ Mac Lane's pentagonal equation $(b5)$.

Our last supposition is the following:
\begin{quote}
(IV)\quad $b^\str$ and $b^\rts$ are isomorphisms inverse to each
other.
\end{quote}
With (I) to (IV) we have obtained exactly the assumptions of
associative categories.

In associative categories, $b^\str$ is 2-lifted to $\mj\kon\mj$,
and 1-lifted to $\mj\kon b^\str$. Either of these assumptions, or
analogous assumptions concerning $b^\rts$, could replace our
assumption (III) for defining associative categories. In
assumption (III), however, just $b^\str$ arrows are involved. The
$\psi^{\Box\kon(\Box\kon\Box),i,\Box}$ and
$\psi^{(\Box\kon\Box)\kon\Box,i,\Box}$ arrows are made of $b^\rts$
arrows when $i=1$, and of both kinds, $b^\str$ and $b^\rts$, when
$i=2$.

Since in symmetric monoidal categories we have $c^m$ together with
$b^\str$ and $b^\rts$ arrows, for every shape $M$ we can conclude
that $\kon$ is $M$-lifting in every shape $N$. Though officially
here $M$ and $N$ are shapes without $\top$, this is easily
extended to shapes with $\top$. The equation $(\psi\alpha*)$ will
hold because we have coherence for symmetric monoidal categories.

To pass from associative categories to monoidal categories, we
need moreover the special object $\top$, the natural isomorphisms
$\delta^\str$, $\delta^\rts$, $\sigma^\str$ and $\sigma^\rts$, and
the equation $(b\delta\sigma)$ of Section~2. (The equation (VII)
can now be derived; see \cite{K64}.) The introduction of the
natural isomorphisms $\delta^\str$, $\delta^\rts$, $\sigma^\str$
and $\sigma^\rts$ may then be justified by suppositions
(IV$\delta$) and (IV$\sigma$), as in Section~2. The equation
$(b\delta\sigma)$ could be justified by understanding the arrows
$\mj_A\kon\sigma_C^\str$ and $\delta_A^\str\kon\mj_C$ in the
diagram
\begin{center}
\begin{picture}(220,80)(0,10)
\put(25,10){$A\kon C$}

\put(0,73){$A\kon(\top\kon C)$}

\put(161,10){$A\kon C$}

\put(156,73){$(A\kon\top)\kon C$}

\put(92,4){$\mj_{A\kon C}$}

\put(60,15){\vector(1,0){90}}

\put(91,86){$b^\str_{A,\top,C}$}

\put(60,78){\vector(1,0){90}}

\put(-2,45){$\mj_A\kon\sigma^\str_C$}

\put(38,65){\vector(0,-1){40}}

\put(178,45){$\delta^\str_A\kon\mj_C$}

\put(173,65){\vector(0,-1){40}}
\end{picture}
\end{center}
\vspace{1ex} as a kind of $\psi$ arrows, going beyond our
$\psi^{N,i,M}$ arrows. The shape $M$ in this extension of $\psi$
arrows would be a \koc\ filled by the ``conjunction of no
formulae'' $\top$. The shapes
$[((\koc\kon\koc)\kon\koc)_\top^2]^L$ and
$[(\koc\kon(\koc\kon\koc))_\top^2]^R$ would both be
$\koc\kon\koc$.

In the presence of $\top$ we can also explain our shapes of the
split $[N_M^i]^S$ in another manner. In monoidal categories the
object $[N_M^i]^L(\vec{B}_{i-1},\vec{A}_m)$ is isomorphic to
$N_M^i(\vec{B}_{i-1},\vec{A}_m,\top,\ldots,\top)$, while
$[N_M^i]^R(\vec{A}_m,\vec{C}_{n-i})$ is isomorphic to
$N_M^i(\top,\ldots,\top,\vec{A}_m,\vec{C}_{n-i})$. In these
isomorphisms only $\delta^\str$, $\delta^\rts$, $\sigma^\str$ and
$\sigma^\rts$ arrows occur essentially. With the extended $\psi$
arrows mentioned above we have that $[N_M^i]^L(\vec{B}_{i-1})$ is
isomorphic to $N_M^i(\vec{B}_{i-1},\top,\top,\ldots,\top)$, the
first $\top$ replacing $\vec{A}_m$, while
$[N_M^i]^R(\vec{C}_{n-i})$ is isomorphic to
$N_M^i(\top,\ldots,\top,\top,\vec{C}_{n-i})$, the last $\top$
replacing $\vec{A}_m$.

\section{\large\bf Coherence for medial commutativity}
The question arises what equations should be assumed for $c^m$ in
the absence of $\top$ in order to obtain coherence in the sense of
Mac Lane's coherence for symmetric monoidal categories. Besides
the categorial equations, the bifunctorial equations for $\kon$,
the naturality and isomorphism equations for $c^m$ and $(\psi
c^m)$, we need equations that would deliver, for example, the
following:
\begin{tabbing}
\mbox{\hspace{2.5em}}\=$(c^m~{\mbox{\rm
I}})$\mbox{\hspace{2em}}\=$c^m\cirk(c^m\kon\mj)\cirk
c^m\cirk(\mj\kon c^m)=(\mj\kon c^m)\cirk c^m\cirk(c^m\kon\mj)\cirk
c^m$
\end{tabbing}
with appropriate indices for the $c^m$ terms so that both sides
are of the type
\[
((A_1\kon A_1')\kon A_2)\kon((A_3\kon A_3')\kon(A_4\kon
A_4'))\str((A_1\kon A_3)\kon A_2)\kon((A_1'\kon A_4)\kon(A_3'\kon
A_4'),
\]
or the following:
\begin{tabbing}
\mbox{\hspace{2.5em}}\=$(c^m~{\mbox{\rm I}})$\mbox{\hspace{2em}}\=
\kill \>$(c^m~{\mbox{\rm II}})$\>$c^m\cirk(c^m\kon\mj)\cirk
c^m\cirk(c^m\kon\mj)\cirk c^m\cirk(c^m\kon\mj)=$\\*[.5ex]
\`$(c^m\kon\mj)\cirk c^m\cirk(c^m\kon\mj)\cirk
c^m\cirk(c^m\kon\mj)\cirk c^m$
\end{tabbing}
with appropriate indices for the $c^m$ terms so that both sides
are of the type
\[
((A_1\kon A_1')\kon(A_2\kon A_2'))\kon((A_3\kon A_3')\kon
A_4)\str((A_1\kon A_1')\kon(A_2\kon A_3'))\kon((A_3\kon A_2')\kon
A_4).
\]
An analogue of $(c^m~{\mbox{\rm II}})$ is obtained by replacing
$c^m\kon\mj$ by $\mj\kon c^m$.

To describe these equations we introduce a syntactically
constructed category we call \Cm. The objects of \Cm\ are the
\emph{formulae} of the propositional language generated by a
nonempty set of \emph{letters} with $\kon$ as sole connective.
Formally, we have the inductive definition:
\begin{quote} every letter is a formula;\\[.5ex]
if $A$ and $B$ are formulae, then $(A\kon B)$ is a formula.
\end{quote}
As usual, we take the outermost parentheses of formulae for
granted, and omit them. We also omit the outermost parentheses of
the \emph{arrow terms} of \Cm, which we define inductively as
follows:
\begin{tabbing}
\mbox{\hspace{2.5em}}\=\mbox{\hspace{2em}}\=\kill

\>for all formulae $A$, $B$, $C$ and $D$,\\[1ex]
\>\>$\mj_A\!:A\str A$ and
$c^m_{A,B,C,D}\!:(A\kon B)\kon(C\kon D)\str(A\kon C)\kon(B\kon
D)$\\[1ex]
\>are arrow terms; if $f\!:A\str B$ and $g\!:C\str D$ are arrow
terms, then\\[1ex]
\>\>$g\cirk f\!:A\str D$, provided $B$ is $C$, and $f\kon
g\!:A\kon C\str B\kon D$\\[1ex]
\>are arrow terms.
\end{tabbing}

These arrow terms are subject to the following equations: the
categorial equations, the bifunctorial equations for $\kon$, the
naturality and isomorphism equations for $c^m$, and, furthermore,
equations to be described in this section. This means that to
obtain \Cm\ we factor the arrow terms through an equivalence
relation engendered by the equations, congruent with respect to
$\kon$ and $\cirk$, and take the equivalence classes as arrows. (A
formal definition of such syntactically constructed categories may
be found in \cite{DP04}, Chapter~2.)

We want to prove that \Cm\ is coherent, which means that there is
a faithful functor $G$ from \Cm\ to the category whose objects are
finite ordinals, with arrows being permutations. For every object
$A$ of \Cm\ we have that $GA$ is the finite ordinal that is the
number of letters in $A$, and on arrows we have definitional
clauses corresponding to the following pictures:

\begin{center}
\begin{picture}(200,50)
\put(8.5,10){\line(0,1){20}} \put(21.5,10){\line(0,1){20}}

\put(15,7){\makebox(0,0)[t]{$GA$}}
\put(15,34){\makebox(0,0)[b]{$GA$}}
\multiput(12,27)(3,0){3}{\circle*{.3}}
\put(0,20){\makebox(0,0)[r]{$G\mj_A$}}

\put(106.5,10){\line(0,1){20}} \put(118.5,10){\line(0,1){20}}
\put(189.5,10){\line(0,1){20}} \put(201.5,10){\line(0,1){20}}
\put(133,10){\line(5,4){25}} \put(148,10){\line(5,4){25}}
\put(158,10){\line(-5,4){25}} \put(173,10){\line(-5,4){25}}

\put(150,34){\makebox(0,0)[b]{$G((A\;\kon\; B)\;\kon\;(C\;\kon\;
D))$}}

\put(150,7){\makebox(0,0)[t]{$G((A\;\kon\; C)\;\kon\;(B\;\kon\;
D))$}} \put(90,20){\makebox(0,0)[r]{$Gc^m$}}

\multiput(109.5,27)(3,0){3}{\circle*{.3}}

\multiput(192.5,27)(3,0){3}{\circle*{.3}}

\multiput(141.5,27)(3,0){3}{\circle*{.3}}

\multiput(158.5,27)(3,0){3}{\circle*{.3}}

\end{picture}
\end{center}

\vspace{2ex}

\begin{center}
\begin{picture}(100,40)
\multiput(30,10)(0,2){10}{\line(0,1){1}}
\multiput(48,10)(0,2){10}{\line(0,1){1}}
\multiput(59,10)(0,2){10}{\line(0,1){1}}
\multiput(77,10)(0,2){10}{\line(0,1){1}}

\put(10,20){\makebox(0,0)[r]{$G(f\kon g)$}}

\put(50,34){\makebox(0,0)[b]{$G(\;A\;\;\kon\;\;C\;)$}}

\put(50,7){\makebox(0,0)[t]{$G(\;B\;\;\kon\;\;D\;)$}}

\put(39,20){\makebox(0,0){$Gf$}} \put(68,20){\makebox(0,0){$Gg$}}
\end{picture}
\end{center}
The functor $G$ maps $\cirk$ into composition of permutations.
(Analogous functors $G$ are defined more formally in \cite{DP04},
Sections 2.9 and 5.1.)

For $n\geq 0$ and $p$ a letter, let $p^n$ be the formula defined
inductively as follows:
\begin{tabbing}
\mbox{\hspace{13em}}\=$p^{n+1}$\= \kill \> $p^0$\> $=p$\\*[.5ex]
\>$p^{n+1}$\> $=p^n\kon p^n$.
\end{tabbing}
We will define inductively for every $n\geq 0$ an auxiliary
category $\Cm(p^n)$ whose only object is $p^n$.

To define the \emph{arrow terms} of $\Cm(p^n)$ we have the
following inductive clauses:
\begin{tabbing}
\hspace{2em}\=for $n\geq 0$, \quad\=$\mj_{p^n}$ is an arrow term
of $\Cm(p^n)$;\\[1ex]\>for $n\geq 2$,\>
$c^m_{p^{n\mn 2},p^{n\mn 2},p^{n\mn 2},p^{n\mn 2}}$ is an arrow
term of
$\Cm(p^n)$;\\[1ex]\>for $n\geq 1$,\>if $f$ and $g$ are arrow terms of
$\Cm(p^{n-1})$, then $f\kon g$ is an ar-\\[.5ex]\>\>row term of
$\Cm(p^n)$;\\[1ex]
\>for $n\geq 0$,\>if $f$ and $g$ are arrow terms of $\Cm(p^n)$,
then $g\cirk f$ is an arrow\\[.5ex]\>\> term of $\Cm(p^n)$.
\end{tabbing}
It is easy to see that $p^n$ is an object of \Cm, and that the
arrow terms of $\Cm(p^n)$ are the arrow terms of \Cm\ of the type
${p^n\str p^n}$. The categories \Cm\ and $\Cm(p^n)$ will be so
defined that $\Cm(p^n)$ is a full subcategory of \Cm.

Let $h(A)$ be the height of the finite binary tree corresponding
to the formula $A$; we take that ${h(p)=0}$. The category \Cm\
will be so defined that there is a faithful functor $F$ from \Cm\
to the category which is the disjoint union of the categories
$\Cm(p^n)$. On objects $FA$ is $p^{h(A)}$. If every letter occurs
in $A$ at most once, then $p^{h(A)}$ is obtained from $A$ by a
uniquely determined uniform substitution; otherwise, the
substitution need not be uniform. For the arrow term ${f\!:A\str
B}$ of \Cm, the arrow term ${Ff\!:FA\str FA}$ of $\Cm(p^{h(A)})$
(note that $FA$ and $FB$ are equal) is obtained by a substitution
in the indices of $f$ induced by the substitution that leads from
$A$ to $FA$.

We will describe the equations we assume for the arrow terms of
$\Cm(p^n)$. The equations of \Cm\ will be obtained from these so
that $F$ is a faithful functor. The procedure for producing the
equations of \Cm\ is the following. We take an equation of
$\Cm(p^n)$, delete its indices, and look for the most general
indices that could be assigned. This gives an equation of \Cm. For
example, the equation $(c^m~{\mbox{\rm I}})$ with the official
indices, mentioned at the beginning of the section, is obtained
from this equation with indices assigned in such a manner that
both sides are of the type $p^3$. (This means that $A_1$, $A_1'$,
$A_3$, $A_3'$, $A_4$ and $A_4'$ are $p$ and $A_2$ is $p^1$, i.e.\
$p\kon p$.)

It remains now to give the equations for the arrow terms of
$\Cm(p^n)$. These equations will make of $\Cm(p^n)$ a group that
is the internal direct product of $G_0^n,\ldots,G_n^n$, where
$G_i^n$, for $0\leq i\leq n$, is isomorphic to the symmetric group
$S_{n \choose i}$, i.e.\ the group of all permutations of $n
\choose i$ elements. The composition $\cirk$ of $\Cm(p^n)$ is
group multiplication, and the identity arrow $\mj_{p^n}$ is the
unit element of the group.

Our notation comes from taking the elements of $G_i^n$ to
correspond to permutations of the set $\{i_1,\ldots,i_{n \choose
i}\}$. The group $G_i^n$ may be generated by the elements that
correspond to the ${n \choose i}\mn 1$ transpositions of immediate
neighbours:
\[
i_{j,j+1}^n\quad\mbox{\rm for}\quad 1\leq j\leq {n \choose i}\mn 1
\]
(with $i_{j,j+1}^n$, the element $i_j$ is mapped to $i_{j+1}$ and
vice versa, while the other elements are fixed). In terms of these
transpositions one defines as usual the transpositions $i_{j,k}^n$
of $G_i^n$ for every $j,k\in\{1,\ldots,{n \choose i}\}$, where
$j<k$:
\begin{tabbing}
\mbox{\hspace{2.5em}}\=$(\emph{def}\,~i_{j,k})$\mbox{\hspace{2em}}\=$i_{j,k}^n=_{df}
i_{j,j+1}^n\cirk \ldots\cirk i_{k-2,k-1}^n\cirk i_{k-1,k}^n\cirk
i_{k-2,k-1}\cirk\ldots\cirk i_{j,j+1}^n.$
\end{tabbing}

We define $\Cm(p^n)$ gradually as follows. (In this definition we
refer to equations given in Section~1.)

\vspace{1ex}

\noindent (0)\quad For the category $\Cm(p^0)$ we assume just the
equation $(\mbox{\it cat}~1)$ for $A$ and $B$ being $p$. This
delivers the categorial equations $(\mbox{\it cat}~1)$ and
$(\mbox{\it cat}~2)$.

\vspace{1ex}

\noindent (1)\quad For the category $\Cm(p^1)$ we assume the
categorial and the bifunctorial equations for $\kon$. As a matter
of fact, the assumption of $(\mbox{\it cat}~2)$ and $(\mbox{\it
bif}~2)$ is here superfluous.

\vspace{1ex}

\noindent (2)\quad For the category $\Cm(p^2)$ we assume the
categorial equations, the bifunctorial equations for $\kon$, and
the naturality and isomorphism equations for $c^m$. As before,
some of these equations are here superfluous---in particular, the
equation $(c^m~\mbox{\it nat})$, which will also be superfluous
for $\Cm(p^n)$ where $n\geq 3$, but which must hold in \Cm. (This
equation is needed when we do not deal with the freely generated
category \Cm, but with categories of that type that may have
additional structure.)

\vspace{1ex}

For $n\geq 2$ we introduce inductively the following
abbreviations:

\[
1_{1,2}^2=_{df} c^m_{p,p,p,p},
\]
\begin{tabbing}
\hspace{2.5em}$(\emph{def}\,~i^{n+1})$\hspace{2em}for $1\leq i\leq
n$,
\end{tabbing}
\vspace{-2ex}
\[
\mbox{\hspace{2em}}i_{j,j+1}^{n+1}=\left\{\begin{array}{ll}
i_{j,j+1}^n\kon\mj_{p^n} & \mbox{\rm if }1\leq j\leq{n \choose
i}\mn 1
\\[1.5ex]
\mj_{p^n}\kon(i\mn 1)_{j-{n \choose i},j-{n \choose i}+1}^n &
\mbox{\rm if }{n \choose i}\pl 1\leq j\leq{{n+1} \choose i}\mn 1.
\end{array}
\right.
\]

We also have the following auxiliary abbreviations:
\[
\begin{array}{rl}
\gamma & \!\!\!\!=_{df} c^m_{p^{n\mn 1},p^{n\mn 1},p^{n\mn
1},p^{n\mn 1}},
\\[1.5ex]
\alpha & \!\!\!\!=_{df} i_{{{n\mn 1} \choose i},{{n\mn 1} \choose
i}+1}^{n+1},
\\[1.5ex]
\beta & \!\!\!\!=_{df} i_{{{n\mn 1} \choose i}+1,{n \choose
i}}^{n+1}.
\end{array}
\]
The right-hand sides of the definitions of $\alpha$ and $\beta$
reduce to abbreviations introduced by $(\emph{def}\,~i^{n+1})$;
with $\beta$, we rely on $(\emph{def}\,~i_{j,k})$ above. Then we
have the following abbreviations:

\begin{tabbing}
\mbox{\hspace{1.8em}}\=({\it norm})
\mbox{\hspace{2em}}\=$i_{k,k+1}^{n+1}\cirk
j_{l,l+1}^{n+1}=j_{l,l+1}^{n+1}\cirk i_{k,k+1}^{n+1}$\kill

\>for $n\geq 3$ and $2\leq i\leq n\mn 1$,\\[1ex]
\>\>$i_{{n \choose i},{n \choose i}+1}^{n+1}=_{df}
\beta\cirk\alpha\cirk\gamma\cirk\alpha\cirk\gamma\cirk\alpha\cirk\beta,$\\[2ex]
\>for $n\geq 2$ and $i=1$,\\[1ex]
\>\>$i_{{n \choose i},{n \choose i}+1}^{n+1}=1_{n,n+1}^{n+1}=_{df}
\alpha\cirk\gamma\cirk\alpha\cirk\gamma\cirk\alpha,$\\[2ex]
\>for $n\geq 2$ and $i=n$,\\[1ex]
\>\>$i_{{n \choose i},{n \choose
i}+1}^{n+1}=n_{1,2}^{n+1}=_{df}n_{2,3}^{n+1}\cirk\gamma\cirk
n_{2,3}^{n+1}\cirk\gamma\cirk n_{2,3}^{n+1}.$
\end{tabbing}

For $n\geq 2$ we have the following.

\vspace{1ex}

\noindent ($n\pl 1$)\quad For the category $\Cm(p^{n+1})$ we first
assume all the equations we have assumed for $\Cm(p^2)$. Next we
assume the equation
\begin{tabbing}
\mbox{\hspace{1.8em}}\=({\it norm})
\mbox{\hspace{2em}}\=$i_{k,k+1}^{n+1}\cirk
j_{l,l+1}^{n+1}=j_{l,l+1}^{n+1}\cirk i_{k,k+1}^{n+1}$, \quad for
$i\neq j$.
\end{tabbing}
As a matter of fact, it is superfluous to assume some instances of
this equation, like, for example,
\[
1_{1,2}^3\cirk 2_{2,3}^3=2_{2,3}^3\cirk 1_{1,2}^3,
\]
which follows from $(\mbox{\it bif}~2)$. Next we assume the
Yang-Baxter equation:
\begin{tabbing}
\mbox{\hspace{1.8em}}\=({\it norm}) \mbox{\hspace{2em}}\=\kill \>
(YB)\> $i_{j,j+1}^{n+1}\cirk i_{j+1,j+2}^{n+1}\cirk
i_{j,j+1}^{n+1}=i_{j+1,j+2}^{n+1}\cirk i_{j,j+1}^{n+1}\cirk
i_{j+1,j+2}^{n+1}$.
\end{tabbing}
With $\prod_{1\leq i\leq n}f(i)$ standing for
$f(1)\cirk\ldots\cirk f(n)$, we assume the equation
\begin{tabbing}
\mbox{\hspace{1.8em}}\=({\it norm}) \mbox{\hspace{2em}}\=\kill
\>$(\emph{def}\,~c^m)$\>$c^m_{p^{n\mn 1},p^{n\mn 1},p^{n\mn
1},p^{n\mn 1}}= \prod_{1\leq i\leq n}(i_{{{n\mn 1} \choose i}+1,{n
\choose i}+1}^{n+1}\cirk \ldots\cirk i_{{n \choose i},{n \choose
i}+{{n\mn 1} \choose {i\mn 1}}}^{n+1})$.
\end{tabbing}
(For the transpositions on the right-hand side we rely again on
$(\emph{def}\,~i_{j,k})$.)

We assume that if $f=g$ holds for $\Cm(p^n)$, then $\mj_{p^n}\kon
f=\mj_{p^n}\kon g$ and ${f\kon\mj_{p^n}=g\kon\mj_{p^n}}$ hold for
$\Cm(p^{n+1})$. As a matter of fact, we need to assume this rule
only when $f=g$ is an instance of $(\emph{def}\,~c^m)$, or is
obtained ultimately from an instance of $(\emph{def}\,~c^m)$
through repeated applications of this rule.

Finally, for $n\geq 3$, we assume also the following equation:
\begin{tabbing}
\mbox{\hspace{1.8em}}\=({\it norm}) \mbox{\hspace{2em}}\=\kill \>
({\it perm})\> $i_{j,j+1}^{n+1}\cirk
i_{j+k,j+k+1}^{n+1}=i_{j+k,j+k+1}^{n+1}\cirk i_{j,j+1}^{n+1}$,
\quad for $k\geq 2$.
\end{tabbing}

\vspace{1.5ex}

This concludes the definition of the categories $\Cm(p^n)$. To
take the easiest nontrivial example, in $\Cm(p^3)$ we find the
following equations. Besides the categorial, bifunctorial,
naturality and isomorphism equations, we have the equation
$(c^m~{\mbox{\rm I}})$ with appropriate indices, which is obtained
by using essentially ({\it norm}). Next we have the equation
$(c^m~{\mbox{\rm II}})$ with appropriate indices, and its analogue
with $c^m\kon\mj$ replaced by $\mj\kon c^m$, which are obtained by
using essentially (YB). Finally, we have the equation $(\psi
c^m)$, which is obtained by using essentially $(\mbox{\it
def}~c^m)$. Conversely, these instances of the equations
$(c^m~{\mbox{\rm I}})$, $(c^m~{\mbox{\rm II}})$, the analogue of
$(c^m~{\mbox{\rm II}})$ and $(\psi c^m)$, together with
categorial, bifunctorial, naturality and isomorphism equations,
are sufficient to define $\Cm(p^3)$. Equations such as
$(c^m~{\mbox{\rm I}})$, $(c^m~{\mbox{\rm II}})$ and $(\psi c^m)$
tend to be rather involved for $\Cm(p^n)$ where $n>3$.

It remains now to verify that, by taking that the composition
$\cirk$ of $\Cm(p^n)$ is group multiplication and that the
identity arrow $\mj_{p^n}$ is the unit element of the group, the
category $\Cm(p^n)$ is a group that is the internal direct product
of $G_0^n,\ldots,G_n^n$, where $G_i^n$, for $0\leq i\leq n$, is
isomorphic to the symmetric group~$S_{n \choose i}$.

It is clear that $\Cm(p^n)$ is a group. Next we specify the
subgroups $G_0^n,\ldots,G_n^n$ of $\Cm(p^n)$ for every $n\geq 0$.
We will always have that $G_0^n$ and $G_n^n$ contain only the
arrow $\mj_{p^n}$, and are the trivial one-element group,
isomorphic to $S_1$, while $G_i^n$ for ${1\leq i\leq n\mn 1}$ is
generated by the arrows $i_{1,2}^n,\ldots,i_{{n \choose i}-1,{n
\choose i}}$, which exist only when $n\geq 2$. By relying on the
bifunctorial equations, and the equations $(\emph{def}\,~i^{n+1})$
and $(\emph{def}\,~c^m)$, we can represent every arrow of
$\Cm(p^n)$ in terms of the abbreviations $i_{j,j+1}^n$, for
${1\leq i\leq n\mn 1}$, the identity arrow $\mj_{p^n}$ and
composition, without $\kon$ and $c^m$.

We arrange the elements of the set $P^n$, which is defined as
\[
\{0_1\}\cup\ldots\cup\{i_1,\ldots,i_{n \choose
i}\}\cup\ldots\cup\{n_1\},
\]
in a sequence $N^n$, defined inductively as follows:
\begin{itemize}
\item[]$N^0$ is $0_1$; \item[]$N^{n+1}$ is the sequence $N^n$
followed by the sequence obtained by replacing $i_j$ in the
sequence $N^n$ by $(i\pl 1)_{j+{n\choose{i+1}}}$ (note that
${n\choose k}=0$ if $k>n$).
\end{itemize}
So we have
\begin{tabbing}
\mbox{\hspace{2.5em}}\=({\it norm}) \mbox{\hspace{2em}}\=\kill
\>\>$N^1=0_11_1$,
\\*[1ex]
\>\>$N^2=0_11_11_22_1$,
\\*[1ex]
\>\>$N^3=0_11_11_22_11_32_22_33_1$,
\\*[1ex]
\>\>$N^4=0_11_11_22_11_32_22_33_11_42_42_53_22_63_33_44_1$, etc.
\end{tabbing}
($N^3$ should be compared with the graphical presentation of
$(\psi c^m)$ in Section~3).

For every arrow ${f\!:p^n\str p^n}$ of $\Cm(p^n)$, and for $G$
being the functor defined earlier in this section, if we replace
$Gp^n$ by $N^n$, then $Gf$ stands for a permutation of $P^n$.
Although somewhat lengthy, it is then straightforward to verify
that $Gi_{j,j+1}^n$ is a transposition of immediate neighbours in
$\{i_1,\ldots,i_{n \choose i}\}$, and that if ${f=g}$ in
$\Cm(p^n)$, then ${Gf=Gg}$.

From that, by relying on the completeness of (YB), ({\it perm})
and equations that hold in virtue of $(c^mc^m)$ for presenting
symmetric groups (see \cite{B11}, Note C, \cite{CM57}, Section
6.2, or \cite{DP04}, Section 5.2), we conclude that $G_i^n$ is
isomorphic to the symmetric group $S_{n \choose i}$, and that
$G_i^n\cap G_j^n$ for $i\neq j$ is $\{\mj_{p^n}\}$. Since we have
the equation ({\it norm}), it follows that $\Cm(p^n)$ is the
internal direct product of $G_0^n,\ldots,G_n^n$.

This implies that ${f=g}$ in $\Cm(p^n)$ if and only if ${Gf=Gg}$.
We have already established the direction from left to right, and
for the other direction we rely on the elimination of $\kon$ and
$c^m$ mentioned above.

We may then prove that the category \Cm\ is coherent, i.e., that
$G$ is a faithful functor. Suppose ${f,g\!:A\str B}$ are arrows of
\Cm\ such that ${Gf=Gg}$. It follows that ${GFf=GFg}$ for $F$
being the faithful functor from \Cm\ to the disjoint union of the
categories $\Cm(p^n)$. The arrows $Ff$ and $Fg$ are elements of
the group $\Cm(p^n)$ for some $n$, and ${GFf=GFg}$ guarantees that
${Ff=Fg}$. From ${Ff=Fg}$ it follows by the faithfulness of $F$
that ${f=g}$.

It is remarkable that when the unit object is added to the
category \Cm, the burden of all the complicated equations of \Cm\
is carried by the equation $(\psi c^m)$.


\begin{thebibliography}{99}

\bibitem{BFSV03} {\sc C. Balteanu, Z. Fiedorowicz, R. Schw\" anzl} and
{\sc R. Vogt}, {\it Iterated monoidal categories}, \textbf
{\textit {Advances in Mathematics}}, vol.\ 176 (2003), pp.\
277-349

\bibitem{B37} {\sc R. Brauer}, {\it On algebras which are connected with the
semisimple continuous groups}, \textbf {\textit {Annals of
Mathematics}}, vol.\ 38 (1937), pp.\ 857-872

\bibitem{B11} {\sc W. Burnside}, \textbf {\textit
{Theory of Groups of Finite Order}}, second edition, Cambridge
University Press, Cambridge, 1911 (reprint, Dover, New York, 1955)

\bibitem{CM57} {\sc H.S.M. Coxeter} and {\sc W.O.J. Moser}, \textbf {\textit
{Generators and Relations for Discrete Groups}}, Springer, Berlin,
1957

\bibitem{DP96} {\sc K. Do\v sen} and {\sc Z. Petri\' c},
{\it Modal functional completeness}, \textbf {\textit {Proof
Theory of Modal Logic}} (H. Wansing, editor), Kluwer, Dordrecht,
1996, pp.\ 167-211

\bibitem{DP03a} --------,
{\it Generality of proofs and its Brauerian representation},
\textbf {\textit {The Journal of Symbolic Logic}}, vol.\ 68
(2003), pp.\ 740-750 (available at: http://
arXiv.\-org/\-math.\-LO/\-0211090)

\bibitem{DP03b} --------, {\it A Brauerian
representation of split preorders}, \textbf {\textit {Mathematical
Logic Quarterly}}, vol.\ 49 (2003), pp.\ 579-586 (available at:
http://\-arXiv.\-org/ math.\-LO/\-0211277)

\bibitem{DP04} --------, \textbf {\textit
{Proof-Theoretical Coherence}}, KCL Publications (College
Publications), London, 2004 (revised version available at:
http://\-www.\-mi.\-sanu. ac.\-yu/\-$\sim$kosta/\-coh.\-pdf)

\bibitem{DP06} --------, {\it Associativity as commutativity},
\textbf {\textit {The Journal of Symbolic Logic}}, vol.\ 71
(2006), pp.\ 217-226 (available at:
http://\-arXiv.\-org/\-math.\-CT/ 0506600)

\bibitem{DP05} --------, {\it Symmetric self-adjunctions:
A justification of Brauer's representation of Brauer's algebras},
\textbf {\textit {Proceedings of the Conference ``Contemporary
Geometry and Related Topics"}} (N. Bokan et al. editors), Faculty
of Mathematics, Belgrade, 2006, pp. 177-187 (available at:
http://\-arXiv.\ org/\-math.\-RT/\-0512102)

\bibitem{EH62} {\sc B. Eckmann} and {\sc P.J. Hilton},
{\it Group-like structures in general categories I:
Multiplications and comultiplications}, \textbf {\textit
{Mathematische Annalen}}, vol.\ 145 (1962), pp.\ 227-255

\bibitem{EK66} {\sc S. Eilenberg} and {\sc G.M. Kelly},
{\it Closed categories}, \textbf {\textit {Proceedings of the
Conference on Categorical Algebra, La Jolla 1965}} (S. Eilenberg
et al., editors), Springer, Berlin, 1966, pp.\ 421-562

\bibitem{JK83} {\sc J. Je\v zek} and {\sc T. Kepka}, \textbf {\textit
{Medial Groupoids}}, Roz\-pra\-vy {\v C}esko\-slovensk{\' e}
Aka\-de\-mie V{\v e}d, {\v R}a\-da ma\-te\-ma\-tick{\' y}ch a p{\v
r}i\-rod\-n{\' \i}ch v{\v e}d, Ro{\v c}\-nik 93, Se{\v s}it 2,
Prague, 1983, 93 pp.

\bibitem{JS93} {\sc A. Joyal} and {\sc R. Street},
{\it Braided tensor categories}, \textbf {\textit {Advances in
Mathematics}}, vol.\ 102 (1993), pp.\ 20-78

\bibitem{K64} {\sc G.M. Kelly},
{\it On Mac Lane's conditions for coherence of natural
associativities, commutativities, etc.}, \textbf {\textit {Journal
of Algebra}}, vol.\ 1 (1964), pp.\ 397-402

\bibitem{ML63} {\sc S. Mac Lane}, {\it Natural associativity and
commutativity}, \textbf {\textit {Rice University Studies, Papers
in Mathematics}}, vol.\ 49 (1963), pp.\ 28-46

\bibitem{ML71} --------, \textbf {\textit {Categories for the Working
Mathematician}}, Springer, Berlin, 1971 (expanded second edition,
1998)

\bibitem{T97} {\sc A. Tonks}, {\it Relating the associahedron
and the permutohedron}, \textbf {\textit {Operads: Proceedings of
Renaissance Conferences}} (J.-L. Loday et al., editors),
Contemporary Mathematics, vol.\ 202, American Mathematical
Society, Providence, 1997, pp.\ 33-36

\end{thebibliography}
\end{document}